\documentclass[twocolumn]{autart}

\usepackage{amsmath}
\usepackage{amssymb}
\usepackage{blkarray}

\newtheorem{theorem}{Theorem}
\newtheorem{lemma}{Lemma}
\newtheorem{assumption}{Assumption}
\newtheorem{remark}{Remark}

\usepackage{algorithm}
\usepackage{algpseudocode}

\usepackage{graphicx}
\graphicspath{ {figure/} }
\usepackage{tabularx}

\usepackage{color}
\usepackage{cite}
\usepackage{url}

\begin{document}
%%%%%%%%%%%%%%%%%%%%%%%%%%%%%%%%%%%%%%%%%%%%%%%%%%%%%%%%%%%%%%%%%%%%%%%%%%%%%%%%
\begin{frontmatter}
\title{A Non-Interior-Point Continuation Method for the Optimal Control Problem with Equilibrium Constraints\thanksref{footnoteinfo}}

\thanks[footnoteinfo]{This work was partially supported by JSPS KAKENHI Grant Number JP22H01510. Kangyu Lin was supported by the CSC Scholarship (No. 201906150138). The material in this paper was partially presented at the 61st IEEE Conference on Decision and Control, December 6 -- 9, 2022, Cancun, Mexico. Corresponding author: Kangyu Lin}

\author[kyotoUniv]{Kangyu Lin}\ead{k-lin@sys.i.kyoto-u.ac.jp},
\author[kyotoUniv]{Toshiyuki Ohtsuka}\ead{ohtsuka@i.kyoto-u.ac.jp}

\address[kyotoUniv]{Department of Systems Science, Graduate School of Informatics, Kyoto University, Kyoto, Japan}

\date{}


\begin{abstract}
    In this study, we focus on the numerical solution method for the optimal control problem with equilibrium constraints (OCPEC).
    It is extremely challenging to solve OCPEC owing to the absence of constraint regularity and strictly feasible interior points.
    To solve OCPEC efficiently, we first relax the discretized OCPEC to recover the constraint regularity and then map its Karush--Kuhn--Tucker (KKT) conditions into a perturbed system of equations.
    Subsequently, we propose a novel two-stage solution method, called the non-interior-point continuation method, to solve the perturbed system.
    In the first stage, a non-interior-point method, which solves the perturbed system using the Newton method and globalizes convergence using a dedicated merit function, is employed.
    In the second stage, a predictor-corrector continuation method is utilized to track the solution trajectory as a function of the perturbed parameter, starting at the solution obtained in the first stage.
    The proposed method regularizes the KKT matrix and does not enforce iterates to remain in the feasible interior, which mitigates the numerical difficulties of solving OCPEC. 
    Convergence properties are analyzed under certain assumptions.
    Numerical experiments demonstrate that the proposed method can accurately track the solution trajectory while demanding significantly less computation time compared to the interior-point method.
    
\end{abstract}

\end{frontmatter}

%%%%%%%%%%%%%%%%%%%%%%%%%%%%%%%%%%%%%%%%%%%%%%%%%%%%%%%%%%%%%%%%%%%%%%%%%%%%%%%%
\section{Introduction}\label{section: introduction}
\subsection{Background and motivation}

Many applied control problems, such as the trajectory optimization of a mechanical system with frictional contacts \cite{Posa2014}, optimal control of a hybrid dynamical system \cite{mynttinen2015smoothing}, and differential Nash equilibrium game \cite{chen2014differential}, consist of dynamical systems that involve equilibria in the form of inequalities (e.g., unilateral constraints) and disjunctive conditions (e.g., frictional contact and mode switching). 
This results in nonsmoothness and discontinuity in the solution trajectories of the dynamical system and introduces challenges in system simulation and control. 
One applicable and general modeling paradigm for this class of nonsmooth dynamical systems, called differential variational inequalities (DVI), is presented in \cite{pang2008differential}. 
DVI consists of an ordinary differential equation (ODE) and a dynamic version of finite-dimensional variational inequalities (VI) \cite{facchinei2003finite}.
VI is a key feature of DVI, which serves as a unifying framework for modeling many types of equilibria. 
A comprehensive study of DVI and VI can be found in \cite{pang2008differential, facchinei2003finite, acary2008numerical}.

Optimal control problems (OCP) aim to search for a state and control pair constrained by a controlled system and path constraints while optimizing a dedicated cost function for some goals.
The most practical approach to solving OCP is the direct method (`first discretize then optimize'), where OCP is first discretized into a finite-dimensional nonlinear programming (NLP) problem, then solved by off-the-shelf NLP solvers.
Based on the types of optimization variables, the direct method can be divided into shooting methods (optimizing state and control) \cite{betts2010practical}, collocation methods (optimizing coefficient of polynomial)\cite{von1993numerical}, and control parametrization methods (optimizing coefficient of control function) \cite{zhu2022sequential, teo2021applied}, all of which demonstrate good control performance.
OCP for a class of nonsmooth systems governed by DVI is referred to as the optimal control problem with equilibrium constraints (OCPEC).
Recently, studies on the optimality conditions of OCPEC have been reported in \cite{Guo2016}, and dedicated numerical solvers based on the direct method have achieved progress in several practical applications \cite{Howell2022, nurkanovic2023time}.

Solving OCPEC is extremely challenging.
First, unless nonsmooth points can be specified, the discretized DVI will always exhibit an incorrect sensitivity, regardless of the chosen discretization step \cite{stewart2010optimal}.
Consequently, NLP solvers may be trapped in the spurious local solution near the initial guess.
Second, the discretized OCPEC is the mathematical programming with equilibrium constraints (MPEC), which is a highly difficult NLP because equilibrium constraints violate most standard constraint qualifications (CQ) at any feasible point.
These CQs, such as the linear independence CQ (LICQ) and the Mangasarian--Fromovitz CQ (MFCQ), guarantee the uniqueness and boundedness of the Lagrangian multipliers, and the structural stability of the optimization problem \cite{Scheel2000}.
Hence, MPEC-tailored solution methods need to be considered.
Additionally, a correct sensitivity can be recovered by smoothing (or relaxing) the system and discretizing it with a sufficiently small timestep \cite{nurkanovic2020limits}. 
This indicates that solving the discretized OCPEC with a decreasing smooth (or relaxation) parameter, i.e., the continuation method, has the potential to be used as a practical approach, and many MPEC-tailored methods are indeed based on the continuation method.

Many MPEC-tailored methods have been proposed over the past three decades, and they can be roughly categorized into relaxation methods \cite{hoheisel2013theoretical},  smoothing methods \cite{jiang2000smooth}, penalization methods \cite{leyffer2006interior}, combinatorial methods \cite{kim2020mpec}, and nonsmooth methods \cite{outrata2013nonsmooth}.
Some degenerate NLP solvers, such as stabilized sequential quadratic programming (SQP) methods \cite{gill2017stabilized}, also demonstrate the potential to solve MPEC.
Considering numerous factors, such as problem size, differentiability, and practical implementation, Scholtes' relaxation methods \cite{scholtes2001convergence} with interior-point-based solvers appear to be the best choice for the MPEC arising from nonsmooth OCP, as concluded in a recent survey \cite{nurkanovic2023solving}.
However, this approach faces numerical difficulties when the relaxation parameter is close to zero, that is, the interior-point-based solver may stall or fail as the feasible interior shrinks towards the disjunctive and empty set.
Furthermore, each problem needs to be solved exactly to track the solution trajectory, hindering the computational efficiency.

\subsection{Contributions}
In this study, to solve the OCPEC efficiently, we first relax the discretized OCPEC and map its Karush--Kuhn--Tucker (KKT) system into a perturbed equation system.
Subsequently, we propose a novel two-stage method to efficiently solve the perturbed equation system.
In the first stage, a non-interior-point method \cite{Lin2022} is utilized to drive the iterate into a region close to the optimum of continuous-time OCPEC.
In the second stage, a predictor-corrector continuation method is employed to efficiently track the solution trajectory.
The most notable feature of the proposed method is the treatment of inequality constraints, which constitutes our core idea and contributes to efficiently solving the OCPEC.
The contributions of our method are as follows:
\begin{itemize}
    \item It can handle a large number of inequalities without enforcing all iterates to remain in the feasible interior, enabling it to address large-scale problems without the feasible interior. 
    \item It generates a perturbed system of equations, enabling the utilization of predictor-corrector continuation methods to track the solution trajectory efficiently.
\end{itemize}
Our core idea is inspired by the fact that the KKT conditions are essentially an equilibrium problem, which can be solved using a class of efficient methods called the non-interior-point path-following method \cite{Kanzow1996} that uses the perturbed complementarity function (C-function).
The most popular C-function is the Fisher--Burmeister (FB) function.
Recent studies have employed several variants of the FB function in the convex QP for real-time applications \cite{liao2018regularized, liao2020fbstab}, and the indirect method for OCP with pure state inequality constraints \cite{fabien2016noninterior}.
Both implementations attain competitive results against interior-point-based solvers.
The method proposed in this study is a generalization of these implementations from the convex and standard OCP to the nonconvex and challenging OCP.
The major differences between the proposed method and existing implementations \cite{liao2018regularized, liao2020fbstab, fabien2016noninterior} lie in two aspects: a more suitable globalization routine for the nonconvex problem and a more efficient continuation method for tracking the solution trajectory.
Additionally, we improve our earlier research \cite{Lin2022} with the predictor-corrector continuation method and characterize the convergence properties.

\subsection{Outline}
The remainder of this paper is organized as follows. 
In Section \ref{section: problem formulation}, we introduce a formulation framework that transforms continuous-time OCPEC into a perturbed NLP problem. 
In Section \ref{section: solution method}, the proposed non-interior-point continuation solution method is presented in detail. 
In Section \ref{section: convergence analysis}, we present the convergence analysis. 
The numerical experiments and results are presented in Section \ref{section: numerical simulation}.
Finally, Section \ref{section: conclusion} presents the conclusions of this study.

\subsection{Notation}
Given the Euclidean $n$-dimensional vector space $\mathbb{R}^n$, we denote its nonnegative orthant by $\mathbb{R}^n_{+}$.
Given a variable $x \in \mathbb{R}^n$, we denote its $i$-th element by $x_i$.
Given a differentiable function $f(x): \mathbb{R}^n \rightarrow \mathbb{R}^m$, we denote the Jacobian matrix of $f$ as $\nabla_x f \in \mathbb{R}^{m \times n}$.
Given a differentiable function $f(x): \mathbb{R}^n \rightarrow \mathbb{R}$, we denote the Hessian matrix of $f$ as $\nabla_{xx} f \in \mathbb{R}^{n \times n}$.
Considering two nonnegative infinite sequences of scalars $\{ f^k\}_{k=1}^{\infty}$ and $\{ g^k\}_{k=1}^{\infty}$, we state that $f^k = O(g^k)$ if there exists a constant $\varepsilon > 0$ such that $| f^k | \leq \varepsilon | g^k |$ holds for all sufficiently large $k$, and we state that $f^k = o(g^k)$ if $ \lim\limits_{k \rightarrow \infty} \frac{f^k}{g^k} = 0$.

%%%%%%%%%%%%%%%%%%%%%%%%%%%%%%%%%%%%%%%%%%%%%%%%%%%%%%%%%%%%%%%%%%%%%%%%%%%%%%%%
\section{Problem Formulation}\label{section: problem formulation}

\subsection{Variational inequalities}
Finite-dimensional VI serves as the unified mathematical formalism for a variety of equilibrium problems, see \cite{facchinei2003finite}.
Its definition is provided below:
\begin{defn}[Variational Inequalities]\label{definition: VI}
    Consider a closed convex nonempty set $K \subseteq \mathbb{R}^{n_\lambda}$, and a continuous function $F: \mathbb{R}^{n_\lambda} \rightarrow \mathbb{R}^{n_\lambda}$, the variational inequalities, denoted VI$(K, F)$, is to find a vector $\lambda \in K$ such that,
    \begin{equation}\label{equation: VI definition}
        (\omega - \lambda)^T F(\lambda) \geq 0, \quad \forall \omega \in K.
    \end{equation}    
    We denote the solution set to VI$(K, F)$ as SOL$(K, F)$.
\end{defn}
In this study, we focus on the \textit{box-constrained} VI (BVI), where set $K$ exhibits a box-constrained structure:
\begin{equation}\label{equation: box constraint VI set K}
    K := \{\lambda \in \mathbb{R}^{n_\lambda} | b_l \leq \lambda \leq b_u \},
\end{equation}
with $b_l \in \{\mathbb{R} \cup \{-\infty\} \}^{n_{\lambda}}$, $b_u \in \{\mathbb{R} \cup \{+\infty\} \}^{n_{\lambda}}$, and $b_l < b_u$.
BVI is common and contains many standard equilibrium problems, for example: 
\begin{itemize}
    \item If $b_l = -\infty, b_u = + \infty$, then BVI is reduced to a \textit{system of equations}: $F(\lambda) = 0$;
    \item If $b_l = 0, b_u = + \infty$, then BVI is reduced to a \textit{complementary problem} (CP): $0 \leq \lambda \perp F(\lambda) \geq 0$.
\end{itemize}
BVI exhibits combinatorial nature, which implies a \textit{switch-case} logical condition: 
If $\lambda = b_l$, then $F(\lambda) \geq 0$; 
If $\lambda = b_u$, then $F(\lambda) \leq 0$; 
If $b_l < \lambda < b_u$, then $F(\lambda) = 0$.
The solution set to the BVI can be explicitly represented, i.e., $\lambda \in$ SOL$(K, F)$ if and only if $\lambda$ satisfies a set of $n_{\lambda}$ equalities and $4 n_{\lambda}$ inequalities:
\begin{subequations}\label{equation: BVI solution set explicit expression}
    \begin{align}
        & F(\lambda) - \eta = 0,\\
        & b_l \leq \lambda \leq b_u, \label{equation: BVI solution set explicit expression auxiliary variable}\\
        & (\lambda - b_l) \odot \eta \leq 0, \label{equation: BVI solution set explicit expression bilinear term 1}\\
        & (b_u - \lambda) \odot \eta \geq 0, \label{equation: BVI solution set explicit expression bilinear term 2}
    \end{align}
\end{subequations}
where $\odot$ is the Hadamard product, and $\eta \in \mathbb{R}^{n_\lambda}$ is an auxiliary variable for function $F$.
The region formed by (\ref{equation: BVI solution set explicit expression auxiliary variable}) -- (\ref{equation: BVI solution set explicit expression bilinear term 2}) includes three pieces: the nonnegative part of axis $\lambda = b_l$, the nonpositive part of axis $\lambda = b_u$, and the segment of axis $\eta = 0$ between $\lambda = b_l$ and $\lambda = b_u$.

\subsection{Differential variational inequalities}
DVI is first proposed in \cite{pang2008differential} to uniformly model a class of nonsmooth dynamical systems. 
It comprises an ODE and a dynamical VI:
\begin{subequations}\label{equation: DVI definition}
    \begin{align}
        & \Dot{x} = f(x(t), u(t), \lambda(t)),  \label{equation: DVI definition ODE} \\
        & \lambda(t) \in \text{SOL}(K, F(x(t), u(t), \lambda(t))),  \label{equation: DVI definition VI}
    \end{align}
\end{subequations}
where $x(t): [0, T] \rightarrow \mathbb{R}^{n_x}$ is the differential state, 
$u(t): [0, T] \rightarrow \mathbb{R}^{n_u}$ is the control input, 
$\lambda(t): [0, T] \rightarrow \mathbb{R}^{n_\lambda}$ is the algebraic variable, 
$f: \mathbb{R}^{n_x} \times \mathbb{R}^{n_u} \times \mathbb{R}^{n_\lambda} \rightarrow \mathbb{R}^{n_x}$ is the ODE function,
and $\text{SOL}(K, F)$ is the solution set of a (dynamical) VI defined by a set $K \subseteq \mathbb{R}^{n_\lambda}$ and a function $F: \mathbb{R}^{n_x} \times \mathbb{R}^{n_u} \times \mathbb{R}^{n_\lambda} \rightarrow \mathbb{R}^{n_\lambda}$.
Since $\lambda(t)$ has to follow the time-varying and multi-valued $\text{SOL}(K, F)$, it may be nonsmooth, discontinuous, or even unbounded w.r.t. time. 
Therefore, substituting $\lambda(t)$ into the ODE (\ref{equation: DVI definition ODE}) may induce the nonsmoothness to $x(t)$, resulting in numerical difficulties in the system simulation and optimal control.
Herein, we adopt a stable discretization method proposed in \cite{pang2008differential},  which discretizes the ODE (\ref{equation: DVI definition ODE}) using the implicit Euler method, while enforcing the satisfaction of VI (\ref{equation: DVI definition VI}) at each time point $t_n$:
\begin{subequations}\label{equation: discretized DVI}
    \begin{align}
        & x_n = x_{n-1} + f(x_n, u_n, \lambda_n)\Delta t,  \label{equation: discretized DVI ODE} \\
        & \lambda_n \in \text{SOL}(K, F(x_n, u_n, \lambda_n)),  \label{equation: discretized DVI VI} 
    \end{align}
\end{subequations}
with $n = 1,2,...,N$, where $x_n\in \mathbb{R}^{n_x}$ and $\lambda_n\in \mathbb{R}^{n_\lambda}$ are the value of $x(t)$ and $\lambda(t)$ at the time point $t_n$, respectively, $u_n \in \mathbb{R}^{n_{u}}$ is the piecewise constant approximation of $u(t)$ in the interval $(t_{n-1}, t_n]$, $N$ is the number of stages, and $\Delta t : = T/N$ is the time step.
Since we adopt a time-stepping discretization method that approximates the overall nonsmooth trajectory by smooth polynomials in the integration interval rather than exactly detecting the nonsmooth point, it can only achieve first-order discretization accuracy \cite{acary2008numerical}.

\subsection{Optimal control problems with equilibrium constraints}\label{subsection: optimal control problems with equilibrium constraints}

Mathematically, the continuous-time OCPEC considered in this study has the form
\begin{subequations}\label{equation: continuous OCPEC}
    \begin{align}
        \min_{x(\cdot), u(\cdot),  \lambda(\cdot)} \ &  L_{T}(x(T)) +  \int_0^T L_{S}(x(t), u(t), \lambda(t)) dt \\
        \text{s.t.} \  &  \Dot{x} = f(x(t), u(t), \lambda(t)), \label{equation: continuous OCPEC ODE} \\
                       &  \lambda(t) \in \text{SOL}(K, F(x(t), u(t), \lambda(t))), \label{equation: continuous OCPEC VI} \\
                       &  G(x(t), u(t)) \geq 0,  \\
                       &  C(x(t), u(t)) = 0, 
    \end{align}
\end{subequations}
where $L_T: \mathbb{R}^{n_x} \rightarrow \mathbb{R}$ and $L_S: \mathbb{R}^{n_x} \times \mathbb{R}^{n_u} \times \mathbb{R}^{n_\lambda} \rightarrow \mathbb{R}$ are the terminal and stage cost function, respectively, 
dynamical system (\ref{equation: continuous OCPEC ODE}) (\ref{equation: continuous OCPEC VI}) is a DVI as defined in (\ref{equation: DVI definition}), 
$G: \mathbb{R}^{n_x} \times \mathbb{R}^{n_u}  \rightarrow \mathbb{R}^{n_G}$ and $C: \mathbb{R}^{n_x} \times \mathbb{R}^{n_u}  \rightarrow \mathbb{R}^{n_C}$ are the inequality and equality path constraints, respectively.
In this study, we numerically solve the continuous-time OCPEC using the direct multiple shooting method \cite{betts2010practical} with DVI discretization (\ref{equation: discretized DVI}), where the equilibrium constraint (\ref{equation: discretized DVI VI}) is replaced with its explicit expression (\ref{equation: BVI solution set explicit expression}).
Since inequalities (\ref{equation: BVI solution set explicit expression auxiliary variable}) -- (\ref{equation: BVI solution set explicit expression bilinear term 2}) inherently violate MFCQ, where no feasible point can strictly satisfy these inequalities, and the gradients of active inequalities are always linearly dependent, we further relax the inequalities (\ref{equation: BVI solution set explicit expression bilinear term 1}) and (\ref{equation: BVI solution set explicit expression bilinear term 2}) to recover the constraint regularity. 
The discretized OCPEC we solve is then given by
\begin{subequations}\label{equation: discretized OCPEC}
    \begin{align}
    \min_{\boldsymbol{x}, \boldsymbol{u}, \boldsymbol{\lambda}, \boldsymbol{\eta}} & L_T(x_N) + \sum^{N}_{n=1}L_{S}(x_n, u_n, \lambda_n) \Delta t, \label{equation: discretized OCPEC cost function}\\
     \text{s.t.}     \quad              & x_{n-1} + f(x_n, u_n, \lambda_n) \Delta t - x_n = 0,  \label{equation: discretized OCPEC ODE}\\
                                        & F(x_n, u_n, \lambda_n) - \eta_n = 0, \label{equation: discretized OCPEC auxilary variable}\\
                                        & b_l \leq \lambda_n \leq b_u, \label{equation: discretized OCPEC set K} \\
                                        & (\lambda_n - b_l) \odot \eta_n \leq s I_{n_{\lambda} \times 1}, \label{equation: discretized OCPEC relax bilinear 1}\\
                                        & (b_u - \lambda_n) \odot \eta_n \geq -s I_{n_{\lambda} \times 1}, \label{equation: discretized OCPEC relax bilinear 2} \\
                                        & G(x_n, u_n) \geq 0, \label{equation: discretized OCPEC path constraint G}\\
                                        & C(x_n, u_n) = 0,  \ n = 1, \dots ,N \label{equation: discretized OCPEC path constraint C},
    \end{align}
\end{subequations}
where $s \geq 0$ is a scalar relaxation parameter, and $\boldsymbol{x} = [x^T_1, \cdots, x^T_N]^T$, $\boldsymbol{u}= [u^T_1, \cdots, u^T_N]^T$, $\boldsymbol{\lambda}= [\lambda^T_1, \cdots, \lambda^T_N]^T$, and $\boldsymbol{\eta}= [\eta^T_1, \cdots, \eta^T_N]^T$ are the vector collecting all states, controls, algebraic variables, and auxiliary variables along the horizon, respectively.

\begin{figure}[!tbp]
    \centering
    \includegraphics[width=0.9\linewidth]{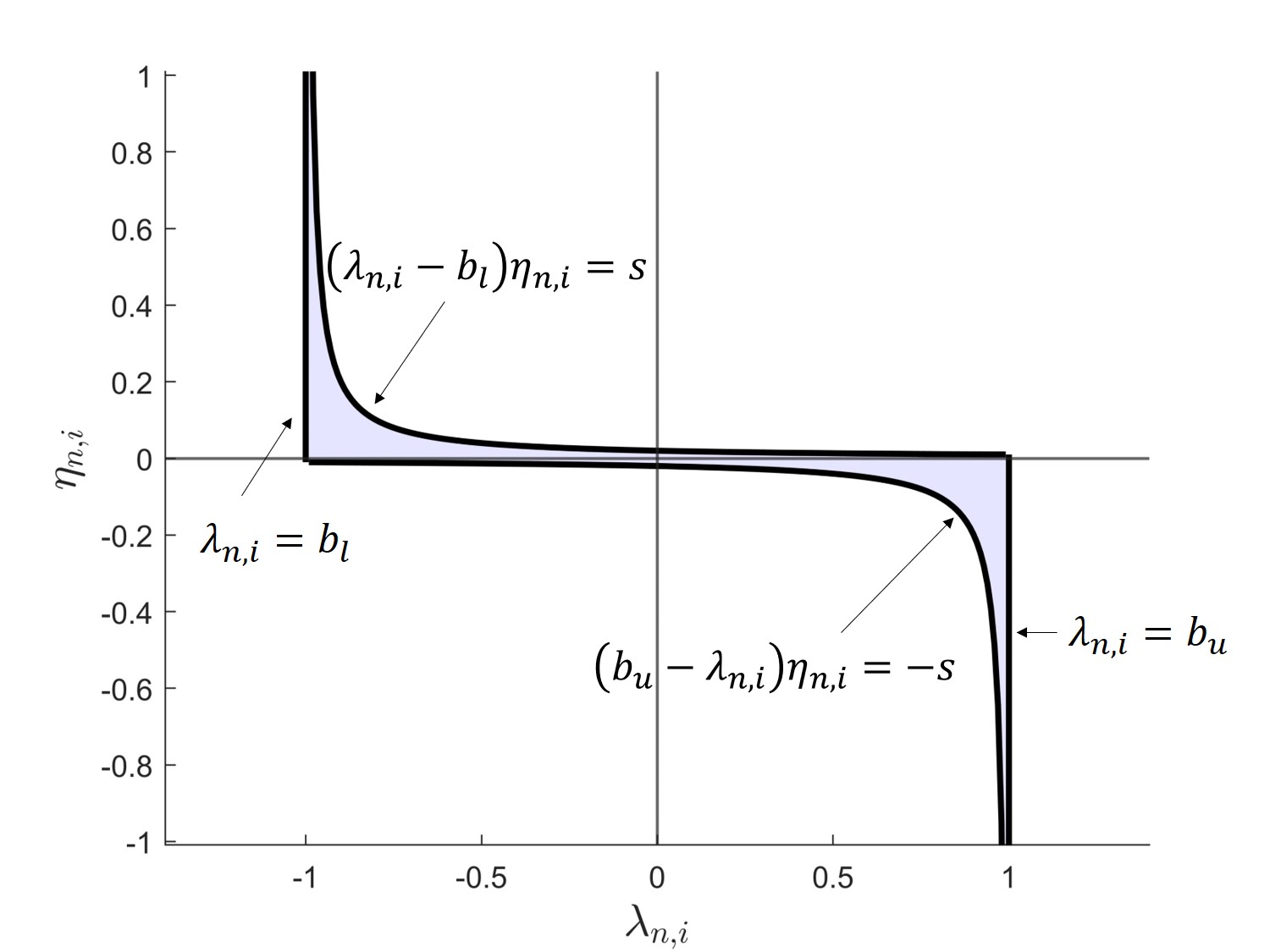}
    \caption{Relaxed feasible set formed by  (\ref{equation: discretized OCPEC set K}) -- (\ref{equation: discretized OCPEC relax bilinear 2}), with $s = 0.02$, $b_l = -1$, and $b_u = 1$. }
    \label{fig: BVI relaxation region}
\end{figure}

As shown in Fig. \ref{fig: BVI relaxation region}, the feasible set formed by relaxed equilibrium constraints (\ref{equation: discretized OCPEC set K}) -- (\ref{equation: discretized OCPEC relax bilinear 2}), exhibits two geometrical features:
First, it always contains the original feasible set regardless of the choice of $s \geq 0$, and reduces to the original set when $s = 0$;
Second, it possesses a strict interior and constraints (\ref{equation: discretized OCPEC set K}) -- (\ref{equation: discretized OCPEC relax bilinear 2}) satisfy the standard CQs when $s > 0$, which enables us to solve the discretized OCPEC (\ref{equation: discretized OCPEC}) with off-the-shelf NLP solvers.
However, two numerical difficulties occur when solving (\ref{equation: discretized OCPEC}).
First, a small time step $\Delta t = o(s)$ needs to be chosen to guarantee the discretization accuracy and the correct sensitivity;
otherwise, the NLP solver may be trapped in certain spurious local solutions far away from the optimum of the continuous-time OCPEC (\ref{equation: continuous OCPEC}).
This is caused by the wrong sensitivity once the condition $\Delta t = o(s)$ (i.e., $\Delta t$ should converge to zero faster than $s$) fails to hold, as discussed in \cite{nurkanovic2020limits} and is further confirmed in Section \ref{subsection: necessity of the continuation method}.
Consequently, a small $\Delta t$ leads to a large-scale problem that renders active-set-based solvers ineffective;
Second, $\Delta t = o(s)$ is a fairly strict condition, and it is more practical to fix $\Delta t$ with a sufficiently small value;
first solve (\ref{equation: discretized OCPEC}) with a relatively large $s \gg \Delta t$ such that the correct sensitivity is recovered and the optimizer can drive iterates into a `good region' close to the optimum of (\ref{equation: continuous OCPEC}), and then perform the continuation method with $s \rightarrow 0$ to refine the iterates.
However, when $s \neq 0$ but $ s \approx 0$, gradients of active inequalities in (\ref{equation: discretized OCPEC set K}) -- (\ref{equation: discretized OCPEC relax bilinear 2}) tend to be linearly dependent. 
Moreover, the feasible interior shrinks towards the disjunctive and empty set, as shown in Fig. \ref{fig: BVI relaxation region}. 
Consequently, interior-point-based solvers may stall or fail to find a new iterate within the feasible interior during the continuation routine.

%%%%%%%%%%%%%%%%%%%%%%%%%%%%%%%%%%%%%%%%%%%%%%%%%%%%%%%%%%%%%%%%%%%%%%%%%%%%%%%%
\section{The Proposed Solution Method}\label{section: solution method}
To mitigate the numerical difficulties mentioned above, we propose a two-stage solution method for solving the discretized OCPEC (\ref{equation: discretized OCPEC}).
The details of the proposed solution method are discussed in the following subsections.

\subsection{Optimality conditions}
For brevity, we collect 
all variables into a vector $\boldsymbol{z} = [z^T_1, \cdots z^T_n, \cdots z^T_N]^T$ with $z_n = [x^T_n, u^T_n, \lambda^T_n, \eta^T_n]^T$,
all equality constraints (\ref{equation: discretized OCPEC ODE}) (\ref{equation: discretized OCPEC auxilary variable}) (\ref{equation: discretized OCPEC path constraint C}) into a vector $\boldsymbol{h} = [h^T_1, \cdots h^T_n, \cdots h^T_N]^T$ with $h_n = [x^T_{n-1} + f^T_n \Delta t - x^T_n, F^T_n - \eta^T_n, C^T_n]^T$, 
inequality constraints (\ref{equation: discretized OCPEC set K}) (\ref{equation: discretized OCPEC path constraint G}) into a vector $\boldsymbol{c} = [c^T_1, \cdots c^T_n, \cdots c^T_N]^T$ with $c_n = [\lambda^T_n - b^T_l, b^T_u - \lambda^T_n, G^T_n]^T$, 
and inequality constraints (\ref{equation: discretized OCPEC relax bilinear 1}) (\ref{equation: discretized OCPEC relax bilinear 2}) into a vector $\boldsymbol{g} = [g^T_1, \cdots g^T_n, \cdots g^T_N]^T$ with $g_n = [s I_{1 \times n_{\lambda}} - (\lambda_n - b_l)^T \odot \eta^T_n, s I_{1 \times n_{\lambda}} + (b_u - \lambda_n)^T \odot \eta^T_n]^T$.
Moreover, we use the shorthand $\boldsymbol{J} = L_T + \sum\limits^{N}_{n=1}L_{S} \Delta t$ for the cost function defined in (\ref{equation: discretized OCPEC cost function}).
The relaxed problem (\ref{equation: discretized OCPEC}) can then be rewritten as a general NLP with a relaxation parameter $s$, denoted as $\mathcal{P}(s)$:
\begin{subequations}\label{equation: general NLP}
    \begin{align}
        \mathcal{P}(s): \quad \min_{\boldsymbol{z}} \ & \boldsymbol{J}(\boldsymbol{z}), \\
        \text{s.t.}     \quad   & \boldsymbol{h}(\boldsymbol{z}) = 0, \\
                                & \boldsymbol{c}(\boldsymbol{z}) \geq 0, \label{equation: general NLP inequalities c}\\
                                & \boldsymbol{g}(\boldsymbol{z}, s) \geq 0 \label{equation: general NLP inequalities g}.
    \end{align}
\end{subequations}
Let $\boldsymbol{\gamma}_h \in \mathbb{R}^{n_h}, \boldsymbol{\gamma}_c \in \mathbb{R}^{n_c}$, and $\boldsymbol{\gamma}_g \in \mathbb{R}^{n_g}$ be the Lagrangian multipliers for the constraints $\boldsymbol{h}$, $\boldsymbol{c}$, and $\boldsymbol{g}$, respectively. 
The Lagrangian of the relaxed problem $\mathcal{P}(s)$ can be defined as:
\begin{equation}
    \mathcal{L}(\boldsymbol{Y}, s) = \boldsymbol{J} + \boldsymbol{\gamma}_h^T\boldsymbol{h} - \boldsymbol{\gamma}^T_c \boldsymbol{c} - \boldsymbol{\gamma}^T_g \boldsymbol{g},
\end{equation}
with $\boldsymbol{Y} = [\boldsymbol{\gamma}^T_h, \boldsymbol{\gamma}^T_c, \boldsymbol{\gamma}^T_g, \boldsymbol{z}^T]^T$ being a vector that collects all primal and dual variables.
The KKT conditions associated with the relaxed problem $\mathcal{P}(s)$ are:
\begin{subequations}\label{equation: KKT condition of relaxed NLP problem}
    \begin{align}
        & \boldsymbol{h} = 0, \\
        & \boldsymbol{\gamma}_c \geq 0, \boldsymbol{c} \geq 0, \boldsymbol{\gamma}_c \odot \boldsymbol{c} = 0, \label{equation: KKT condition of relaxed NLP problem complementary condition c}\\
        & \boldsymbol{\gamma}_g \geq 0, \boldsymbol{g} \geq 0, \boldsymbol{\gamma}_g \odot \boldsymbol{g} = 0, \label{equation: KKT condition of relaxed NLP problem complementary condition g}\\
        & \nabla_{\boldsymbol{z}} \mathcal{L} = 0.
    \end{align}
\end{subequations}
For a given $s > 0$, a point $\boldsymbol{Y}^*(s)$ that satisfies (\ref{equation: KKT condition of relaxed NLP problem}) is referred to as the KKT point of $\mathcal{P}(s)$, and its primal part $\boldsymbol{z}^*(s)$ is referred to as the stationary point of $\mathcal{P}(s)$.
We make the following standard assumptions for $\mathcal{P}(s)$:
\begin{assumption}\label{assumption: relaxed NLP problem}
    For any given $s> 0$, $\mathcal{P}(s)$ satisfies the following conditions:
    \begin{itemize}
        \item (A1) Functions $\boldsymbol{J}, \boldsymbol{h}, \boldsymbol{c}$, and $\boldsymbol{g}$ are all twice Lipschitz continuously differentiable;
        \item (A2) There exists at least one KKT point $\boldsymbol{Y}^*(s)$;
        \item (A3) The Jacobian of equality constraints $\nabla_{\boldsymbol{z}}\boldsymbol{h}$ has a full row rank at the KKT point $\boldsymbol{Y}^*(s)$,
        \item (A4) The reduced Hessian $W^T \nabla_{\boldsymbol{z} \boldsymbol{z}} \mathcal{L} W $ is positive definite at the KKT point $\boldsymbol{Y}^*(s)$, where $W \in \mathbb{R}^{n_z \times (n_z - n_h)}$ has a full column rank and its column is the basis for the null space of $\nabla_{\boldsymbol{z}}\boldsymbol{h}$.
    \end{itemize}
\end{assumption}
Here, (A1) and (A2) guarantee the differentiability and feasibility of $\mathcal{P}(s)$.
(A3) is slightly weaker than LICQ as we do not make the full-rank assumption on the Jacobian of active inequalities in $\nabla_{\boldsymbol{z}}\boldsymbol{c}$ and $\nabla_{\boldsymbol{z}}\boldsymbol{g}$.
(A4) is required to guarantee the nonsingularity of the KKT matrix, as stated in Proposition \ref{proposition: nonsingularity of the KKT matrix}.

\begin{remark}[Stationarity Properties]
    In \cite{hoheisel2013theoretical}, it was proven that under certain assumptions (e.g., MPEC-LICQ), $\lim\limits_{s\rightarrow 0} \boldsymbol{z}^*(s) = \boldsymbol{z}^*$ and the accumulation point $\boldsymbol{z}^*$ is a C-stationary point of $\mathcal{P}(0)$. 
    This constitutes the core idea of the continuation method, which approximately solves $\mathcal{P}(0)$ by solving a sequence of relatively easier NLP problem $\mathcal{P}(s)$ with $s \rightarrow 0$. 
    Note that the C-stationary point is a weaker stationarity notion than the KKT point and may have a descent direction, which is a limitation of the relaxation method. 
    Refer to Theorem 3.1 in \cite{hoheisel2013theoretical} for details.
\end{remark}

\subsection{Treatments of inequality constraints}
As discussed above, $\mathcal{P}(s)$ is a large-scale problem, making it essential to efficiently treat inequality constraints (\ref{equation: general NLP inequalities c}) and (\ref{equation: general NLP inequalities g}).
In this study, we treat inequality constraints using the smoothed FB function \cite{Kanzow1996}.
Specifically, complementary conditions (\ref{equation: KKT condition of relaxed NLP problem complementary condition c}) and (\ref{equation: KKT condition of relaxed NLP problem complementary condition g}) are respectively mapped into the following perturbed systems of equations:
\begin{subequations} \label{equation: KKT condition reformulated complementary condition}
    \begin{align}
        & \Psi_{c} := \Psi(\boldsymbol{\gamma}_c, \boldsymbol{c}, \sigma) = 0, \\
        & \Psi_{g} := \Psi(\boldsymbol{\gamma}_g, \boldsymbol{g}, \sigma) = 0.
    \end{align}
\end{subequations}
Here, the function $\Psi$ is defined in an element-wise manner by the smoothed FB function denoted as $\psi$:
\begin{equation}\label{equation: smoothed FB function}
   \psi (a, b, \sigma) = \sqrt{a^2 + b^2 + \sigma^2}-a-b,
\end{equation}
with scalar variables $a, b$ and a smoothed parameter $\sigma \geq 0$. 
The smoothed FB function $\psi$ has an important property:
\begin{equation}\label{smoothFB_properties}
    \psi (a, b, \sigma) = 0 \; \Leftrightarrow \; a \geq 0,\; b \geq 0, \; ab = \frac{1}{2} \sigma^2.
\end{equation}
Moreover, the function $\psi$ is smooth for any $\sigma > 0$, and nonsmooth only when $a = b = \sigma = 0$.
We collect all parameters into a vector $p = [s, \sigma]^T$.
The KKT system (\ref{equation: KKT condition of relaxed NLP problem}) can then be mapped into a perturbed system of equations with a parameter vector $p$, as stated below:
\begin{equation}\label{equation: nonsymmetric KKT condition}
    \boldsymbol{T}(\boldsymbol{Y}, p) = 
    \begin{bmatrix}
        \boldsymbol{h}^T & \Psi^T_{c} & \Psi^T_{g} & \nabla_{\boldsymbol{z}} \mathcal{L}^T
    \end{bmatrix}^T
    = 0.
\end{equation}
Let $\boldsymbol{Y}_{\sigma}^*(s)$ be the solution to $\boldsymbol{T}(\boldsymbol{Y}, p)$, then for any given $s > 0$, 
$\boldsymbol{Y}_{\sigma}^*(s)$ is equivalent to $\boldsymbol{Y}^*(s)$ when $\sigma = 0$, thereby, we can obtain $\boldsymbol{Y}^*(s)$ by solving a sequence of smooth perturbed equation system (\ref{equation: nonsymmetric KKT condition}) with $\sigma \rightarrow 0$. 
The smoothed FB function was chosen for two reasons, which also form our core idea to mitigate the numerical difficulties of solving $\mathcal{P}(s)$.
First, it treats the inequality constraints efficiently similar to the interior-point methods, i.e., considering all the inequalities at each iteration and following a central path $ab = \frac{1}{2} \sigma^2$.
However, it enforces $a, b\geq 0$ \textit{implicitly} via the FB function value rather than \textit{explicitly} via the feasible initial guess and fraction-to-the-boundary rule.
This allows iterates to be infeasible, resulting in a larger stepsize and search space against geometrical difficulties of the feasible set.
Second, the KKT system (\ref{equation: KKT condition of relaxed NLP problem}) is mapped into a perturbed system of equations, which can be solved efficiently by using Newton's method with sparsity structure exploitation.
Moreover, the solution trajectory $\boldsymbol{Y}^*(s)$ can be tracked efficiently by tracking $\boldsymbol{Y}_{\sigma}^*(s)$ using the predictor-corrector continuation method rather than solving each relaxed problem $\mathcal{P}(s)$ \textit{exactly} as in the case of the interior-point methods.

In the following subsections, we introduce details of the proposed two-stage method for solving the challenging problem $\mathcal{P}(s)$ with a sufficiently small $s^J$. 
The proposed solution method is summarized in Algorithm \ref{algorithm: non-interior-point continuation method}.

\begin{algorithm}[!tbp]
\caption{Non-interior-point continuation method to solve the relaxed problem $\mathcal{P}(s)$}
\label{algorithm: non-interior-point continuation method}
\begin{algorithmic}
    \State \textbf{Input: } $\{ p^j = [s^j, \sigma^j]^T \}^J_{j = 0}$, $\boldsymbol{Y}^0$.
    \State \textbf{Output: } $\boldsymbol{Y}^*$. 
\end{algorithmic}

\begin{algorithmic}[1]
    \State \textbf{Stage 1: solve $\mathcal{P}(s)$ with $p^0$}.
    \State $\boldsymbol{Y}^{*, 0} \gets $ \text{NON-INTERIOR-POINT} $(\boldsymbol{Y}^0, p^0)$ 
    \State \textbf{Stage 2: track $\boldsymbol{Y}^*(s)$ along $\{ p^j \}^J_{j = 0}$}.
    \State $\boldsymbol{Y}^{*, J} \gets $ \text{CONTINUATION} $(\{p^{j}\}^{J}_{j=0}, \boldsymbol{Y}^{*, 0})$
    \State $\boldsymbol{Y}^* \gets \boldsymbol{Y}^{*, J}$
\end{algorithmic}

\end{algorithm}

\subsection{Stage 1: Non-interior-point solution method}\label{subsection: first stage method}

In the first stage, we aim to drive iterates into a `good region' close to the optimum of the continuous-time OCPEC (\ref{equation: continuous OCPEC}) by solving a relatively easier problem $\mathcal{P}(s)$ with a larger parameter $s^0 > 0$.
As discussed in subsection \ref{subsection: optimal control problems with equilibrium constraints}, this requires $\mathcal{P}(s)$ to exhibit a correct sensitivity, which can be achieved by specifying a small $\Delta t$ and a large $s^0$ such that $\Delta t \ll s^0$.
In the first stage, we do not update but fix the parameter vector $p$ as $p^0 = [s^0, \sigma^0]^T$ with a relatively large smoothed parameter $\sigma^0 > 0$ for improved numerical performance.
We then utilize a non-interior-point method \cite{Lin2022} to efficiently solve the perturbed equation system (\ref{equation: nonsymmetric KKT condition}) with a fixed $p^0$, which is summarized in Algorithm \ref{algorithm: non-interior-point method} and detailed below.

\begin{algorithm}[!tbp]
\caption{Non-interior-point method}
\label{algorithm: non-interior-point method}
\begin{algorithmic}
    \Function{non-interior-point}{$\boldsymbol{Y}^0, p^0$}
    \State \textbf{Initialization}: $\boldsymbol{Y}^k \gets \boldsymbol{Y}^0$.
    \For{$k = 1 \text{ to } k_{max}$}
    \State \textbf{Step 1. } \text{check termination condition: }
    \If{ (\ref{equation: termination condition}) is satisfied}
        \State $\boldsymbol{Y}^{*} \gets \boldsymbol{Y}^k$.
        \State break iteration routine.
    \EndIf
    \State \textbf{Step 2. } \text{evaluate KKT residual and matrix: }
    \State $\boldsymbol{T}^k, \mathcal{K}^k \gets \boldsymbol{Y}^k, p^0$ \Comment{(\ref{equation: nonsymmetric KKT condition}) (\ref{equation: nonsymmetric KKT matrix})}
    \State \textbf{Step 3. } \text{evaluate search direction: }
    \State $\Delta \boldsymbol{Y} \gets \boldsymbol{T}^k, \mathcal{K}^k$. \Comment{(\ref{equation: nonsymmetric linear equation system})}
    \State \textbf{Step 4. } \text{merit line search: }
    \State $\boldsymbol{Y}^{k + 1} \gets \boldsymbol{Y}^k, \Delta \boldsymbol{Y}, p^0$. \Comment{(\ref{equation: Armijo condition})}
    \EndFor
    \State \Return $\boldsymbol{Y}^{*}$
    \EndFunction
\end{algorithmic}

\end{algorithm}

Let the current iterate be $\boldsymbol{Y}^{k}$, we define the primal residual $E^k_p$, dual residual $E^k_d$, and complementary residual $E^k_c$ associated with $\boldsymbol{Y}^{k}$, respectively, as:
\begin{subequations}\label{equation: optimality error}
    \begin{align}
        & E^k_p = \max (\| \boldsymbol{h}^k\|_{\infty}, \| \min(0,  [(\boldsymbol{c}^k)^T, (\boldsymbol{g}^k)^T]) \|_{\infty}), \\
        & E^k_d = \frac{1}{\kappa_d}\max (\|\nabla_{\boldsymbol{z}} \mathcal{L}^k\|_{\infty}, \| \min(0, [(\boldsymbol{\gamma}^k_c)^T, (\boldsymbol{\gamma}^k_g)^T]) \|_{\infty}), \\
        & E^k_c = \frac{1}{\kappa_c} \max(\| \boldsymbol{c}^k \odot \boldsymbol{\gamma}^k_c  \|_{\infty}, \| \boldsymbol{g}^k \odot \boldsymbol{\gamma}^k_g  \|_{\infty}) ,
    \end{align}
\end{subequations}
where $\kappa_d$, $\kappa_c \geq 1$ are scaling parameters defined in \cite{wachter2006implementation}.
$\boldsymbol{Y}^{k}$ is returned as the optimal solution of the first stage if one of the following termination conditions are satisfied:
\begin{subequations}\label{equation: termination condition}
    \begin{align}
        & E^k_{kkt} := \max (E^k_p, E^k_d,  E^k_c) \leq \epsilon_{kkt}, \label{equation: termination condition composite KKT}\\
        &  E^k_p \leq \epsilon_p, \ E^k_d \leq \epsilon_d, \ E^k_c \leq \epsilon_c, \label{equation: termination condition individual KKT}
    \end{align}
\end{subequations}
where $\epsilon_p, \epsilon_d, \epsilon_c, \epsilon_{kkt} > 0$ are the given tolerances.
(\ref{equation: termination condition composite KKT}) is the default setting.
(\ref{equation: termination condition individual KKT}) is introduced because $E^k_d$  and  $E^k_c$ may fail to satisfy a composite tolerance $ \epsilon_{kkt}$, as $E^k_c$ is controlled by $\sigma$, as shown in (\ref{smoothFB_properties}), and some elements in $\boldsymbol{\gamma}^k_c, \boldsymbol{\gamma}^k_g$ may be extremely large, owing to the nearly linear dependence of the gradients of almost active inequalities.

If the termination condition (\ref{equation: termination condition}) is not satisfied, then a new search direction $\Delta \boldsymbol{Y}$ is evaluated by applying Newton's method to solve the KKT system (\ref{equation: nonsymmetric KKT condition}).
This leads to the following linear equation system:
\begin{equation}\label{equation: nonsymmetric linear equation system}
    \mathcal{K}(\boldsymbol{Y}^k, p^0) \Delta \boldsymbol{Y} = - \boldsymbol{T}(\boldsymbol{Y}^k, p^0),
\end{equation}
where the KKT matrix $\mathcal{K}(\boldsymbol{Y}, p)$ has the  form:
\begin{equation}\label{equation: nonsymmetric KKT matrix}
    \begin{split}
        & \mathcal{K}(\boldsymbol{Y}, p) := \nabla_{\boldsymbol{Y}} \boldsymbol{T}(\boldsymbol{Y}, p) = \\
        &
        \begin{bmatrix}
        -\nu_h I & 0 & 0 & \nabla_{\boldsymbol{z}}\boldsymbol{h} \\
        0  & \nabla_{\boldsymbol{\gamma}_c} \Psi_{c} -\nu_c I & 0 & \nabla_{\boldsymbol{c}} \Psi_{c} \nabla_{\boldsymbol{z}} \boldsymbol{c} \\
        0  &  0 & \nabla_{\boldsymbol{\gamma}_g} \Psi_{g} -\nu_g I& \nabla_{\boldsymbol{g}} \Psi_{g} \nabla_{\boldsymbol{z}} \boldsymbol{g} \\
        \nabla_{\boldsymbol{z}}\boldsymbol{h}^T & -\nabla_{\boldsymbol{z}}\boldsymbol{c}^T & -\nabla_{\boldsymbol{z}}\boldsymbol{g}^T & \nabla_{\boldsymbol{z} \boldsymbol{z}} \mathcal{L} + \nu_H I\\        
        \end{bmatrix}\\
    \end{split},    
\end{equation}
with $\nabla_{\boldsymbol{\gamma}_c} \Psi_{c}$, $\nabla_{\boldsymbol{c}} \Psi_{c}$, $\nabla_{\boldsymbol{\gamma}_g} \Psi_{g}$, and $\nabla_{\boldsymbol{g}} \Psi_{g}$ being the diagonal matrices. 
Moreover, they are all negative definite when $\sigma > 0$, as $\nabla_a \psi, \nabla_b \psi < 0$ holds when $\sigma > 0$ (see Lemma \ref{lemma: smoothed FB function properties}).
Here, $\nu_h$, $\nu_c$, $\nu_g$, and $\nu_H \geq 0$ are the regularization parameters that regularize the singularity (if Assumption \ref{assumption: relaxed NLP problem} fails to hold, see Proposition \ref{proposition: nonsingularity of the KKT matrix}) or the ill-condition of the KKT matrix $\mathcal{K}$ for better numerical performance.

Upon obtaining a new search direction $\Delta \boldsymbol{Y}$, we evaluate the stepsize $\alpha$ by the merit line search procedure inspired by \cite{Waltz2006} to guarantee the convergence from a poor initial guess.
A dedicated merit function $\Theta(\boldsymbol{Y},p^0)$ is defined to simultaneously measure the cost and constraint violation, which is a $\ell_1$ non-differentiable exact penalty function expressed below:
\begin{subequations}\label{equation: merit function}
    \begin{align}
        & \Theta(\boldsymbol{Y},p^0) = \boldsymbol{J}(\boldsymbol{z}) + \beta \| \boldsymbol{M}(\boldsymbol{Y},p^0) \|_1, \\
        & \boldsymbol{M}(\boldsymbol{Y},p^0) = 
        \begin{bmatrix}
            \boldsymbol{h}^T & \Psi_{c}^T & \Psi_{g}^T
        \end{bmatrix}^T,
    \end{align}
\end{subequations}
where $\boldsymbol{M}(\boldsymbol{Y},p^0)$ represents the constraint violation, and $\beta > 0$ is a penalty parameter whose choice is stated below.
Let $\mathbb{D}\Theta^{k}$ denote the directional derivative of the merit function (\ref{equation: merit function}) at the current iterate $\boldsymbol{Y}^{k}$ along the given direction $\Delta \boldsymbol{Y}$. 
To guarantee the existence of the stepsize $\alpha$ to be uniformly bounded away from zero (see Theorem \ref{theorem: existence of non zero stepsize}) and provide a sufficient decrease for the merit function, we enforce $\mathbb{D}\Theta^{k}$ to be sufficiently negative:
\begin{equation}\label{equation: directional derivative of merit function}
    \mathbb{D}\Theta^{k} = \nabla_{\boldsymbol{z}} \boldsymbol{J}^{k} \Delta \boldsymbol{z} - \beta \| \boldsymbol{M}^{k}\|_1 \leq - \rho \beta \| \boldsymbol{M}^{k}\|_1,
\end{equation}
with a parameter $\rho \in (0, 1)$, where the equality in (\ref{equation: directional derivative of merit function}) is the explicit expression of $\mathbb{D}\Theta^{k}$, owing to the properties of $\ell_1$ function (see Theorem 18.2 in \cite{nocedal2006numerical}).
We can then choose a penalty parameter $\beta > 0$ satisfying
\begin{equation}\label{equation: penalty parameter in merit function}
    \beta \geq \beta_{trail} := \frac{\nabla_{\boldsymbol{z}} \boldsymbol{J}^{k} \Delta \boldsymbol{z}}{(1 - \rho) \|\boldsymbol{M}^{k}\|_1},
\end{equation}
and perform a backtracking line search strategy to evaluate a new iterate $\boldsymbol{Y}^{k + 1} = \boldsymbol{Y}^{k} + \alpha \Delta \boldsymbol{Y}$.
The trail stepsize $\alpha \leftarrow \nu_{\alpha} \alpha$ is gradually reduced from $\alpha_{max} = 1$ with $\nu_{\alpha} \in (0, 1)$ until the Armijo condition is satisfied
\begin{equation}\label{equation: Armijo condition}
    \Theta(\boldsymbol{Y}^{k} + \alpha \Delta \boldsymbol{Y}, p^0) \leq \Theta(\boldsymbol{Y}^k, p^0) + \nu_D \alpha \mathbb{D}\Theta^{k},
\end{equation}
with $\nu_D \in (0, 1)$ representing the desired reduction in $\Theta$.
In our setting $\rho = 0.1$, $ \nu_{\alpha} = 0.5$, and $\nu_D = 10^{-4}$.

\subsection{Stage 2: Predictor-corrector continuation method}
In the second stage, we aim to evaluate a solution for the more difficult problem $\mathcal{P}(s)$ with a smaller $s^J$ by tracking the solution trajectory $\boldsymbol{Y}^*(s)$, starting at a solution of the easier problem $\mathcal{P}(s)$ with a larger $s^0$.
Tracking refers to computing a sequence of points $\{\boldsymbol{Y}^{*,j}\}^{J}_{j=0}$ with $\boldsymbol{Y}^{*,j}$ being an approximation of $\boldsymbol{Y}^*(s)$ at $s^j$. 
Instead of existing methods that track $\boldsymbol{Y}^*(s)$ by solving each problem $\mathcal{P}(s)$ exactly, we track $\boldsymbol{Y}^*(s)$ efficiently by utilizing the predictor-corrector continuation method to track the perturbed solution trajectory $\boldsymbol{Y}_{\sigma}^*(s)$, starting at the solution obtained in the first stage based on $p^0$. 
This includes one predictor step and multiple corrector steps, as summarized in Algorithm \ref{algorithm: predictor-corrector continuation method}.

\begin{algorithm}[!tbp]
\caption{Predictor-corrector continuation method}
\label{algorithm: predictor-corrector continuation method}   
\begin{algorithmic}
\Function{continuation}{$\{ p^j\}^J_{j = 0}$, $\boldsymbol{Y}^{*, 0}$}
    \State \textbf{Initialization: } $\boldsymbol{Y}^{*, j} \gets \boldsymbol{Y}^{*, 0}$
    \For{$j = 0 \text{ to } J - 1$}
    \State \textbf{Step 1. } \text{Euler predictor step: }
    \State $\boldsymbol{Y}^{*, j+1}_{Eu} \gets \boldsymbol{Y}^{*, j}, p^j, p^{j + 1}$. \Comment{(\ref{equation: Euler predictor step})}
    \State \textbf{Step 2. } \text{Newton corrector step: }
    \State $\boldsymbol{Y}^{*, j+1}_{Ne} \gets \boldsymbol{Y}^{*, j+1}_{Eu}, p^{j+1}$. \Comment{(\ref{equation: Newton corrector step})}
    \State \textbf{Step 3. } \text{Additional Newton corrector step: }
    \If{ additional corrector step is required}
    \For{$i = 1 : \nu_{Ne}$}
    \State $\boldsymbol{Y}^{*, j+1}_{Ne} \gets \boldsymbol{Y}^{*, j+1}_{Ne}, p^{j + 1}$. \Comment{(\ref{equation: additional Newton corrector step})}
    \EndFor
    \EndIf
    \State \textbf{Step 4. } \text{prepare for the next iteration: }
    \State $\boldsymbol{Y}^{*, j+1} \gets \boldsymbol{Y}^{*, j+1}_{Ne}$.
    \EndFor 
\State \Return $\boldsymbol{Y}^{*, J}$
\EndFunction
\end{algorithmic}

\end{algorithm}

Given the current iterate $\boldsymbol{Y}^{*,j}$ and parameter vector $p^j$, to evaluate a new iterate $\boldsymbol{Y}^{*, j + 1}$ at the new parameter vector $p^{j + 1}$, we first perform an Euler predictor step to obtain the Euler-iterate $\boldsymbol{Y}^{*, j+1}_{Eu}$:
\begin{subequations}\label{equation: Euler predictor step}
    \begin{align}
        & \Delta p = p^{j + 1} - p^{j},\\
        & \Delta \boldsymbol{Y}_{Eu} =  - \mathcal{K}(\boldsymbol{Y}^{*, j}, p^j)^{-1} \mathcal{S}(\boldsymbol{Y}^{*, j}, p^j) \Delta p,\\
        & \boldsymbol{Y}^{*, j+1}_{Eu} = \boldsymbol{Y}^{*, j} + \Delta \boldsymbol{Y}_{Eu},
    \end{align}
\end{subequations}
where $\mathcal{S}(\boldsymbol{Y}, p):= \nabla_{p} \boldsymbol{T}(\boldsymbol{Y}, p)$ is the parameter sensitivity matrix having the form
\begin{equation}\label{equation: sensitivity matrix}
    \mathcal{S}(\boldsymbol{Y}, p) = 
    \begin{bmatrix}
        0 & 0 \\
        0 & \nabla_{\sigma} \Psi_{c} \\
        \nabla_{\boldsymbol{g}} \Psi_{g} \nabla_{s} \boldsymbol{g} & \nabla_{\sigma} \Psi_{g} \\
        0 & 0
    \end{bmatrix},
\end{equation}
with $\nabla_{s} \boldsymbol{g} = I_{n_{g} \times 1}$.
Based on the computed Euler-iterate $\boldsymbol{Y}^{*, j+1}_{Eu}$, we then perform a Newton corrector step to obtain the Newton-iterate $\boldsymbol{Y}^{*, j+1}_{Ne}$:
\begin{subequations}\label{equation: Newton corrector step}
    \begin{align}
        & \Delta \boldsymbol{Y}_{Ne} = - \mathcal{K}(\boldsymbol{Y}^{*, j+1}_{Eu}, p^{j + 1})^{-1} \boldsymbol{T}(\boldsymbol{Y}^{*, j+1}_{Eu}, p^{j + 1}), \\
        & \boldsymbol{Y}^{*, j+1}_{Ne} = \boldsymbol{Y}^{*, j+1}_{Eu} + \Delta \boldsymbol{Y}_{Ne}.
    \end{align}
\end{subequations}
Furthermore, additional Newton corrector steps can be performed repeatedly to refine the solution accuracy:
\begin{subequations}\label{equation: additional Newton corrector step}
    \begin{align}
        & \Delta \boldsymbol{Y}^{+}_{Ne} = - \mathcal{K}(\boldsymbol{Y}^{*, j+1}_{Ne}, p^{j + 1})^{-1} \boldsymbol{T}(\boldsymbol{Y}^{*, j+1}_{Ne}, p^{j + 1}), \\
        & \boldsymbol{Y}^{*, j+1}_{Ne} \leftarrow \boldsymbol{Y}^{*, j+1}_{Ne} + \Delta \boldsymbol{Y}^{+}_{Ne}.
    \end{align}
\end{subequations}
However, the number of additional Newton corrector steps denoted by $\nu_{Ne}$, should be carefully specified based on the trade-off between efficiency and accuracy.

%%%%%%%%%%%%%%%%%%%%%%%%%%%%%%%%%%%%%%%%%%%%%%%%%%%%%%%%%%%%%%%%%%%%%%%%%%%%%%%%
\section{Convergence Analysis}\label{section: convergence analysis}

\subsection{Nonsingularity of the KKT matrix}
We first investigate the nonsingularity of the KKT matrix $\mathcal{K}$ defined by (\ref{equation: nonsymmetric KKT matrix}).
This is critical because (\ref{equation: nonsymmetric linear equation system}), (\ref{equation: Euler predictor step}), (\ref{equation: Newton corrector step}), and (\ref{equation: additional Newton corrector step}) require computing $\mathcal{K}^{-1}$.
Moreover, it guarantees path regularity when tracking the solution trajectory of a perturbed equation system \cite{allgower2012numerical}.
\begin{theorem}\label{proposition: nonsingularity of the KKT matrix}
    Let Assumption \ref{assumption: relaxed NLP problem} hold and pick $\sigma > 0$.
    Then, for any given $s \geq 0$, the KKT matrix $\mathcal{K}(\boldsymbol{Y}, p)$ with $p = [s, \sigma]^T$ is nonsingular for any $\boldsymbol{Y} \in \mathbb{R}^{n_Y}$ in the neighborhood of a KKT point $\boldsymbol{Y}^*(s)$.
\end{theorem}
\begin{pf}
    Since the Assumption \ref{assumption: relaxed NLP problem} holds, we can set $\nu_h$, $\nu_c$, $\nu_g$, and $\nu_H = 0$.
    Suppose the KKT matrix $\mathcal{K}$ is singular, then there exists a non-zero vector $q \in \mathbb{R}^{n_Y}$ such that $\mathcal{K} q = 0$.
    By dividing $q = [q^T_1, q^T_2, q^T_3, q^T_4]^T$ with $q_1 \in \mathbb{R}^{n_h}$, $q_2 \in \mathbb{R}^{n_c}$, $q_3 \in \mathbb{R}^{n_g}$, and $q_4 \in \mathbb{R}^{n_z}$, we obtain:
    \begin{subequations}\label{equation: KKT matrix times vector}
        \begin{align}
            & \nabla_{\boldsymbol{z}}\boldsymbol{h} q_4 = 0, \label{equation: KKT matrix times vector 1}\\
            & \nabla_{\boldsymbol{\gamma}_c} \Psi_{c} q_2 + \nabla_{\boldsymbol{c}} \Psi_{c} \nabla_{\boldsymbol{z}} \boldsymbol{c} q_4 = 0, \label{equation: KKT matrix times vector 2}\\
            & \nabla_{\boldsymbol{\gamma}_g} \Psi_{g} q_3 + \nabla_{\boldsymbol{g}} \Psi_{g} \nabla_{\boldsymbol{z}} \boldsymbol{g} q_4 = 0, \label{equation: KKT matrix times vector 3}\\
            & \nabla_{\boldsymbol{z}}\boldsymbol{h}^T q_1 -\nabla_{\boldsymbol{z}}\boldsymbol{c}^T q_2 -\nabla_{\boldsymbol{z}}\boldsymbol{g}^T q_3 + \nabla_{\boldsymbol{z} \boldsymbol{z}} \mathcal{L} q_4 = 0 \label{equation: KKT matrix times vector 4}.
        \end{align}
    \end{subequations}
    Substituting (\ref{equation: KKT matrix times vector 2}) and (\ref{equation: KKT matrix times vector 3}) into (\ref{equation: KKT matrix times vector 4}) to eliminate $q_2$, and $q_3$, (\ref{equation: KKT matrix times vector}) becomes
    \begin{equation}\label{equation: simplified KKT matrix times vector}
        \begin{bmatrix}
            \nabla_{\boldsymbol{z}}\boldsymbol{h}^T & \nabla_{\boldsymbol{z} \boldsymbol{z}} \mathcal{L} + \mathcal{R}_c + \mathcal{R}_g \\
                0 & \nabla_{\boldsymbol{z}}\boldsymbol{h}
        \end{bmatrix}
        \begin{bmatrix}
            q_1 \\
            q_4
        \end{bmatrix}
        = 0,        
    \end{equation}
    with matrices $\mathcal{R}_c, \mathcal{R}_g$ being
    \begin{subequations}
        \begin{align*}
            & \mathcal{R}_c = \nabla_{\boldsymbol{z}}\boldsymbol{c}^T (\nabla_{\boldsymbol{\gamma}_c} \Psi_{c}^{-1} \nabla_{\boldsymbol{c}} \Psi_{c}) \nabla_{\boldsymbol{z}}\boldsymbol{c}, \\
            & \mathcal{R}_g = \nabla_{\boldsymbol{z}}\boldsymbol{g}^T (\nabla_{\boldsymbol{\gamma}_g} \Psi_{g}^{-1} \nabla_{\boldsymbol{g}} \Psi_{g})  \nabla_{\boldsymbol{z}}\boldsymbol{g}.
        \end{align*}
    \end{subequations}
    Since $\nabla_{\boldsymbol{\gamma}_c} \Psi_{c}$, $\nabla_{\boldsymbol{c}} \Psi_{c}$, $\nabla_{\boldsymbol{\gamma}_g} \Psi_{g}$, and $\nabla_{\boldsymbol{g}} \Psi_{g}$ are all negative definite diagonal matrices when $\sigma > 0$ (see Lemma \ref{lemma: smoothed FB function properties}),
    we have $\mathcal{R}_c, \mathcal{R}_g \succ 0$ when $\sigma > 0$.
    By left-multiplying both sides in (\ref{equation: simplified KKT matrix times vector}) with vector $[q_4^T,q_1^T]$, we obtain:
    \begin{equation} \label{equation: hessian times q 4}
        q^T_4 (\nabla_{\boldsymbol{z} \boldsymbol{z}} \mathcal{L} + \mathcal{R}_c + \mathcal{R}_g) q_4 = 0.
    \end{equation}
    Since $q_4$ lies in the null space of $\nabla_{\boldsymbol{z}}\boldsymbol{h}$, there exists a vector $\xi \in \mathbb{R}^{n_z - n_h}$ such that $q_4 = W \xi$, and (\ref{equation: hessian times q 4}) becomes
    \begin{equation}
        \xi^T W^T \nabla_{\boldsymbol{z} \boldsymbol{z}} \mathcal{L} W \xi^T + \xi^T W^T(\mathcal{R}_c + \mathcal{R}_g) W \xi = 0.
    \end{equation}
    Since $W^T \nabla_{\boldsymbol{z} \boldsymbol{z}} \mathcal{L} W, \mathcal{R}_c, \mathcal{R}_g \succ 0$, we have $\xi = 0$, and accordingly $q_4 = 0$.
    Hence, we have $q_1 = 0$ following from (\ref{equation: simplified KKT matrix times vector}) as  $\nabla_{\boldsymbol{z}}\boldsymbol{h}$ being full row rank, and $q_2 = q_3 = 0$ following from (\ref{equation: KKT matrix times vector 2}) (\ref{equation: KKT matrix times vector 3}), which contradicts the assumption made at the beginning that $q$ is a non-zero vector. \qed
\end{pf}

\subsection{Global convergence}
In this subsection, we prove the existence of a non-zero stepsize $\alpha$ in the merit line search procedure.

\begin{theorem}\label{theorem: existence of non zero stepsize}
    Given an iterate $\boldsymbol{Y}^k$ and a new search direction $\Delta \boldsymbol{Y}$ evaluated by (\ref{equation: nonsymmetric linear equation system}) based on $\boldsymbol{Y}^k$.
    Suppose the constraint violation $\boldsymbol{M}(\boldsymbol{Y}^k, p^0) \neq 0$.
    Then, there always exists a stepsize $\alpha \in (0, 1]$ that satisfies the Armijo condition (\ref{equation: Armijo condition}).
\end{theorem}
\begin{pf}
    We temporarily ignore the argument $p$ of the merit function $\Theta$ and constraint violation $\boldsymbol{M}$ to streamline the proof.
    By subtracting $\alpha \mathbb{D}\Theta^{k}$ from both sides of the inequality (\ref{equation: Armijo condition}), we obtain:
    \begin{equation}\label{equation: equivalent formulation of Armijo condition}
        \Theta(\boldsymbol{Y}^k + \alpha \Delta \boldsymbol{Y}) - \Theta(\boldsymbol{Y}^k) - \alpha \mathbb{D}\Theta^{k} \leq (\nu_D - 1) \alpha \mathbb{D}\Theta^{k}.
    \end{equation}
    Regarding the term $\Theta(\boldsymbol{Y}^k + \alpha \Delta \boldsymbol{Y}) - \Theta(\boldsymbol{Y}^k)$, based on Taylor's theorem and triangle inequalities, we have:
    \begin{equation}\label{equation: upper bound for the difference of merit function}
        \begin{split}
            & \Theta(\boldsymbol{Y}^k + \alpha \Delta \boldsymbol{Y}) - \Theta(\boldsymbol{Y}^k) \\
            = \ & \nabla_{\boldsymbol{z}} \boldsymbol{J}^{k} \alpha \Delta \boldsymbol{z} + \int^1_0 (\nabla_{\boldsymbol{z}} \boldsymbol{J}(\boldsymbol{z}^k + t \alpha \Delta \boldsymbol{z}) - \nabla_{\boldsymbol{z}} \boldsymbol{J}^k) \alpha \Delta \boldsymbol{z} dt \\
             & - \beta \| \boldsymbol{M}^k \|_1 + \beta \|  \boldsymbol{M}^k + \nabla_{\boldsymbol{Y}} \boldsymbol{M}^k \alpha \Delta \boldsymbol{Y} \\
             & + \int^1_0  (\nabla_{\boldsymbol{Y}} \boldsymbol{M}(\boldsymbol{Y}^k + t \alpha \Delta \boldsymbol{Y}) - \nabla_{\boldsymbol{Y}} \boldsymbol{M}^k) \alpha \Delta \boldsymbol{Y} dt \|_1 \\
            \leq \ & \nabla_{\boldsymbol{z}} \boldsymbol{J}^{k} \alpha \Delta \boldsymbol{z} - \beta \| \boldsymbol{M}^k \|_1 + \beta \|(1-\alpha) \boldsymbol{M}^k \|_1 \\
            & + \int^1_0 (\nabla_{\boldsymbol{z}} \boldsymbol{J}(\boldsymbol{z}^k + t \alpha \Delta \boldsymbol{z}) - \nabla_{\boldsymbol{z}} \boldsymbol{J}^k) \alpha \Delta \boldsymbol{z} dt \\
            & + \| \int^1_0  (\nabla_{\boldsymbol{Y}} \boldsymbol{M}(\boldsymbol{Y}^k + t \alpha \Delta \boldsymbol{Y}) - \nabla_{\boldsymbol{Y}} \boldsymbol{M}^k) \alpha \Delta \boldsymbol{Y} dt \|_1 \\
            \leq \ & \alpha \mathbb{D}\Theta^{k} + \frac{1}{2}\alpha^2(L_{J} \|\Delta \boldsymbol{z} \|^2_1 + \beta L_M \| \Delta \boldsymbol{Y}\|^2_1),
        \end{split}
    \end{equation}
    where $L_J, L_M > 0$ are the Lipschitz constant for the Jacobian $\nabla_{\boldsymbol{z}}\boldsymbol{J}$ and $\nabla_{\boldsymbol{Y}}\boldsymbol{M}$ under the Lipschitz continuity assumption (A1), respectively.
    The last inequality in (\ref{equation: upper bound for the difference of merit function}) provides an upper bound for the left side of (\ref{equation: equivalent formulation of Armijo condition}).
    Hence, following from (\ref{equation: equivalent formulation of Armijo condition}) and (\ref{equation: upper bound for the difference of merit function}), we only need to prove the existence of stepsize $\alpha \in (0, 1]$ satisfying
    \begin{equation}\label{equation: upper bound satisfy Armijo condition}
        \frac{1}{2}\alpha^2(L_{J} \|\Delta \boldsymbol{z} \|^2_1 + \beta L_M \| \Delta \boldsymbol{Y}\|^2_1) \leq (\nu_D - 1) \alpha \mathbb{D}\Theta^{k}.
    \end{equation}
    Since $\beta > 0$ and is specified to guarantee $\mathbb{D}\Theta^{k} < 0$ when $\boldsymbol{M}^k \neq 0$, together with $\nu_D \in (0, 1)$ and $L_{J} \|\Delta \boldsymbol{z} \|^2_1 + \beta L_M \| \Delta \boldsymbol{Y}\|^2_1 > 0$, we can find a stepsize $\alpha \in (0, 1]$ such that (\ref{equation: upper bound satisfy Armijo condition}) holds, in particular, any $\alpha \leq \hat{\alpha}$ with $\hat{\alpha}$ being
    \begin{equation}
        \hat{\alpha} = \frac{2(\nu_D - 1) \mathbb{D}\Theta^{k}}{L_{J} \|\Delta \boldsymbol{z} \|^2_1 + \beta L_M \| \Delta \boldsymbol{Y}\|^2_1}. 
    \end{equation}  
    Since we reduce the stepsize $\alpha \leftarrow \nu_{\alpha} \alpha$ from $\alpha_{max} = 1$ with $\nu_{\alpha} \in (0, 1)$, the proof is complete. \qed  
\end{pf}

\subsection{Local convergence}
Since in the first stage, we utilize $\ell_1$ non-differentiable exact penalty function (\ref{equation: merit function}), the stepsize can not be guaranteed to recover to one in the neighborhood of the solution.
However, in the second stage, we perform the full-step Newton corrector steps, which exhibit a quadratic local convergence as stated below.

\begin{theorem}\label{theorem: local convergence}
    Let Assumption \ref{assumption: relaxed NLP problem} hold and pick $\sigma > 0$.
    Let $\boldsymbol{Y}^*_{\sigma^j}(s^j)$ be a point of the solution trajectory $\boldsymbol{Y}_{\sigma}^*(s)$ at $p^j = [s^j, \sigma^j]^T$.
    Then, there exists a neighborhood $\mathcal{Y}(\boldsymbol{Y}^*_{\sigma^j}(s^j))$ of $\boldsymbol{Y}^*_{\sigma^j}(s^j)$ such that for every Euler-iterate $\boldsymbol{Y}^{*, j}_{Eu}$ in $ \mathcal{Y}(\boldsymbol{Y}^*_{\sigma^j}(s^j))$, the iterate sequence generated by the Newton corrector steps (\ref{equation: Newton corrector step}) (\ref{equation: additional Newton corrector step}) remains in $\mathcal{Y}(\boldsymbol{Y}^*_{\sigma^j}(s^j))$, and converges to $\boldsymbol{Y}^*_{\sigma^j}(s^j)$ quadratically. 
\end{theorem}

\begin{pf}
    Theorem \ref{proposition: nonsingularity of the KKT matrix} implies the boundness of $\|\mathcal{K}^{-1}\|$.
    Together with the Lipschitz continuity property of $\mathcal{K}$ guaranteed by (A1) in Assumption \ref{assumption: relaxed NLP problem}, the quadratic convergence can be proved, which is a standard result of Newton's method (see Theorem 11.2 in \cite{nocedal2006numerical}). \qed
\end{pf}

\subsection{Error analysis}
Finally, we investigate the solution error induced by the last parameter $s^J, \sigma^J$ of sequence $\{ p^j\}_{j=0}^{J}$.
First, we provide some properties of the smoothed FB function.
\begin{lemma}\label{lemma: smoothed FB function properties}
    The smoothed FB function defined by (\ref{equation: smoothed FB function}) has the following properties
    \begin{subequations}
        \begin{align}   
            & \nabla_{a}\psi (a, b, \sigma), \nabla_{b}\psi (a, b, \sigma) \in [-2, 0], \label{equation: smoothed FB properties gradient} \\
            & | \psi (a, b_1, \sigma_1) - \psi (a, b_2, \sigma_2) | \leq 2 |b_1 - b_2| +  |\sigma_1 - \sigma_2 |. \label{equation: smoothed FB properties error sigma}
        \end{align}
    \end{subequations}     
    
\end{lemma}
\begin{pf}
     The property (\ref{equation: smoothed FB properties gradient}) follows the expressions of $\nabla_{a}\psi (a, b, \sigma), \nabla_{b}\psi (a, b, \sigma)$, and it is also evident that $\nabla_{a}\psi (a, b, \sigma), \nabla_{b}\psi (a, b, \sigma) \in (-2, 0)$ when $\sigma > 0$.
     The property (\ref{equation: smoothed FB properties error sigma}) follows the triangle inequalities:
    \begin{equation}
        \begin{split}
            & | \psi (a, b_1, \sigma_1) - \psi (a, b_2, \sigma_2) | \\
            = & \ | \sqrt{a^2 + b_1^2 + \sigma_1^2} - \sqrt{a^2 + b_2^2 + \sigma_2^2} - (b_1 - b_2)|\\
            \leq & \ \sqrt{(b_1 - b_2)^2 + (\sigma_1 - \sigma_2)^2} + | b_1 - b_2 | \\
            \leq & \ 2| b_1 - b_2 | + | \sigma_1 - \sigma_2 |. \qed
        \end{split}
    \end{equation}     
\end{pf}

\begin{theorem}\label{theorem: solution error boundness}
    Given $s^J > s^{0_+} > 0$ and $\sigma^J > \sigma^{0_+} > 0$ with $s^{0_+}, \sigma^{0_+}$ the parameters that sufficiently close to zero. 
    Let the KKT residual with $p^{0_+} = [s^{0_+}, \sigma^{0_+}]^T$ be defined as $r^{0_+}_{kkt}(\boldsymbol{Y}) = \| \boldsymbol{T}(\boldsymbol{Y}, p^{0_+}) \|_{\infty}$.
    Let $\boldsymbol{Y}^{*, J}$ be the solution of $\boldsymbol{T}(\boldsymbol{Y}, p^J) = 0$ with $p^J = [s^J, \sigma^J]^T$.    
    Subsequently, $r^{0_+}_{kkt}(\boldsymbol{Y}^{*, J})$ is bounded by the given $s^J, \sigma^J$:
    \begin{equation}\label{equation: solution error boundness}
        r^{0_+}_{kkt}(\boldsymbol{Y}^{*, J}) = O(| s^J - s^{0_+} |) + O(|\sigma^J - \sigma^{0_+}|).
    \end{equation}
\end{theorem}

\begin{pf}
    The fourth term $\nabla_{\boldsymbol{z}} \mathcal{L}$ of $\boldsymbol{T}(\boldsymbol{Y}, p)$ does not depend on $s$ and $\sigma$, because $s$ linearly enters $\boldsymbol{g}(\boldsymbol{z}, s) \geq 0$ only and hence does not exist in $\nabla_{\boldsymbol{z}}\boldsymbol{g}$.
    The residual of the second term $\Psi_{c} = \Psi(\boldsymbol{\gamma}_c, \boldsymbol{c}(\boldsymbol{z}), \sigma)$ is smaller than that of the third term $\Psi_{g} = \Psi(\boldsymbol{\gamma}_g, \boldsymbol{g}(\boldsymbol{z}, s), \sigma) $, as $\Psi_{c}$ depends only on $\sigma$.
    Hence, only the third term $\Psi_{g}$ needs to be checked:
    \begin{equation}\label{equation: solution error boundness proof}
        \begin{split}
                & \| \boldsymbol{T}(\boldsymbol{Y}^{*, J}, p^{0_+}) \|_{\infty} \\
            = \ & \| \boldsymbol{T}(\boldsymbol{Y}^{*, J}, p^{0_+}) - \boldsymbol{T}(\boldsymbol{Y}^{*, J}, p^J) \|_{\infty} \\
            \leq \ & 2 \| \boldsymbol{g}(\boldsymbol{z}^{*, J}, s^{0_+}) - \boldsymbol{g}(\boldsymbol{z}^{*, J}, s^J) \|_{\infty} + |\sigma^J - \sigma^{0_+} | \\
            = \ & 2 | s^J - s^{0_+} | + |\sigma^J - \sigma^{0_+} |,
        \end{split}       
    \end{equation}
    where the inequality in (\ref{equation: solution error boundness proof}) follows from Lemma \ref{lemma: smoothed FB function properties}. \qed
\end{pf}

%%%%%%%%%%%%%%%%%%%%%%%%%%%%%%%%%%%%%%%%%%%%%%%%%%%%%%%%%%%%%%%%%%%%%%%%%%%%%%%%
\section{Numerical Simulation}\label{section: numerical simulation}
In this section, we test the non-interior-point continuation method. 
The proposed method was implemented based on the CasADi symbolic framework \cite{Andersson2019} and MATLAB built-in sparse linear algebra routines.  
Simulations were performed in the MATLAB R2023b environment on a laptop PC with a 1.80 GHz Intel Core i7-8550U. 
The implementation is available at \url{https://github.com/KY-Lin22/NIPOCPEC}. 
For comparison, we adopted a well-developed interior-point filter solver called IPOPT \cite{wachter2006implementation} through the CasADi interface.

\subsection{Benchmark problem and implementation details}
\begin{figure}[!tbp]
    \centering
    \includegraphics[width=0.9\linewidth]{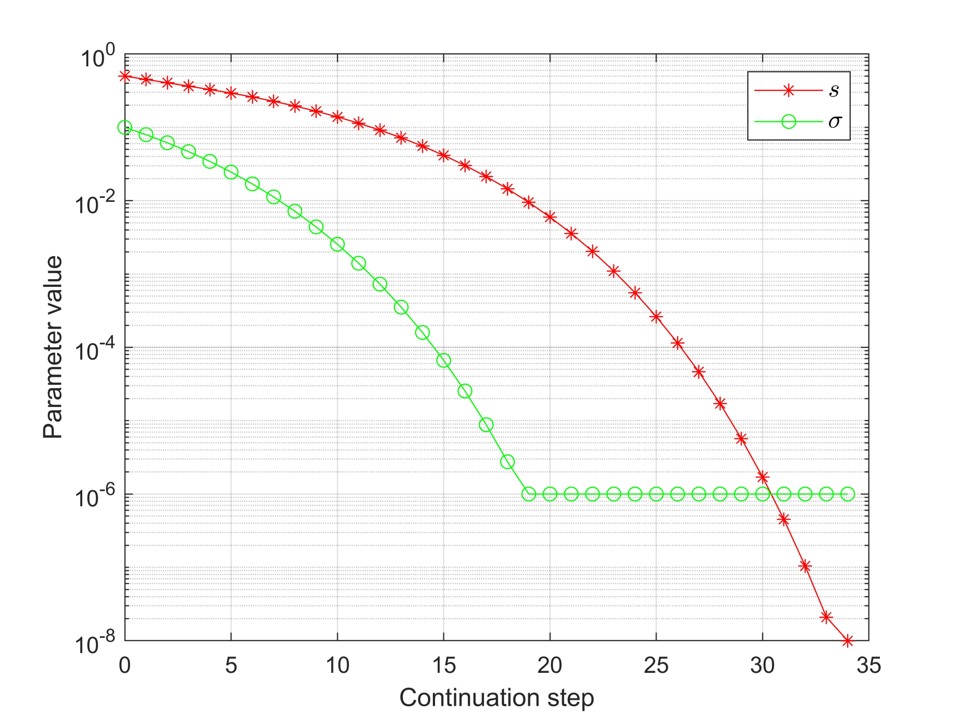}
    \caption{Parameter sequence $\{ p^j = [s^j, \sigma^j]^T\}^J_{j = 0}$}
    \label{fig: parameter sequence}
\end{figure}

We consider OCPEC (\ref{equation: continuous OCPEC}) with a quadratic cost function:
\begin{subequations}
    \begin{align*}
        L_{T} = \ & \| x(T) - x_{e} \|^2_{Q_T},
        \\
        L_{S} = \ & \| x(t) - x_{ref}(t) \|^2_{Q_x} + \|u(t)\|^2_{Q_{u}} + \|\lambda(t) \|^2_{Q_{\lambda}},
    \end{align*}
\end{subequations}
where $Q_T, Q_x, Q_{u}$, and $Q_{\lambda}$ are weighting matrices, $x_{e}$ is the terminal state, and $x_{ref}(t)$ is the reference trajectory.
We consider an example called \textit{Cart Pole with Friction}, which is taken from \cite{Howell2022}.
A cart pole with friction is a nonlinear and nonsmooth dynamic system governed by DVI, which is constructed by two generalized coordinates (cart position $x_c$ and pole angle $\theta_p$), one control input $u_c$ exerting in the cart, and one Coulomb friction $\lambda$ between cart and ground. 
We define the system state as $x = [x_c, \theta_p, \Dot{x}_c, \Dot{\theta}_p]^T$, subsequently, the dynamics are given by
\begin{subequations}\label{equation: affine DVI}
    \begin{align*}
        & \Dot{x} = f(x, u, \lambda) =  
        \begin{bmatrix}
            \Dot{x}_c \\
            \Dot{\theta}_p \\
            \begin{bmatrix}
                m_1 + m_2 & m_2 l \cos(\theta_p) \\
                m_2 l \cos(\theta_p) & m_2 l^2
            \end{bmatrix}^{-1}
            \begin{bmatrix}
                H_1\\
                H_2
            \end{bmatrix}
        \end{bmatrix}, \\
        & H_1 = u + \lambda + m_2 l \sin(\theta_p) (\Dot{\theta}_p)^2, \ H_2 = - m_2 g l \sin(\theta_p),\\
        & F = \Dot{x}_c, \ K = \{\lambda \ |  -2 \leq \lambda \leq 2 \},
    \end{align*}
\end{subequations}
where $m_1 = 1, m_2 = 0.1, l = 1$, and $g = 9.8$.
The goal is to generate a swing motion of the pole from $\theta_p = 0$ to $\theta_p = \pi$ by driving the cart with $u_c$, while the system is subjected to the Coulomb friction $\lambda$. 
Here, we specify $x_{ref}(t) \equiv x_{e}$, $T = 3$, $\nu_h = \nu_c = \nu_g = 10^{-7}$, $\nu_H = 10^{-6}$, $\epsilon_p = 10^{-6}, \epsilon_d = 10^{-4}, \epsilon_c = (\sigma^0)^2$, and $\epsilon_{kkt} = 10^{-6}$.
We specify $s^0 = 5 \cdot 10^{-1}, s^J = 10^{-8}$, $\sigma^0 = 10^{-1}$, and $\sigma^J = 10^{-6}$, and generate the parameter sequence $\{ p^j\}^J_{j = 0}$ by
\begin{subequations}
    \begin{align}
        & s^{j+1} = \max\{s^J, \min \{\kappa_t s^j, (s^j)^{\kappa_e}\} \}, \\
        & \sigma^{j+1} = \max\{\sigma^J, \min \{\kappa_t \sigma^j, (\sigma^j)^{\kappa_e}\} \},
    \end{align}
\end{subequations}
with $\kappa_t = 0.9$, and $\kappa_e = 1.1$, as shown in Fig. \ref{fig: parameter sequence}.
We measure the violation of equilibrium constraints using the \textit{natural residual} $r_{na} = \lambda - \Pi_{K}(\lambda - F)$ with the Euclidean projector $\Pi_{K}$, and $r_{na} = 0$ if and only if $\lambda \in \text{SOL}(K, F)$.
We approximate the Lagrangian Hessian using the Gauss-Newton method, i.e., $\nabla_{\boldsymbol{z} \boldsymbol{z}} \mathcal{L} \approx \nabla_{\boldsymbol{z} \boldsymbol{z}} \boldsymbol{J}$.
We denote the non-interior-point continuation method by \textit{NIP}, and the interior-point continuation method by \textit{IP} (i.e., solve each $\mathcal{P}(s)$ along the sequence $\{s^j\}_{j=0}^{J}$ with the default setting of IPOPT).

\subsection{Solution trajectory}
\begin{figure}[!tbp]
    \centering
    \includegraphics[width=0.9\linewidth]{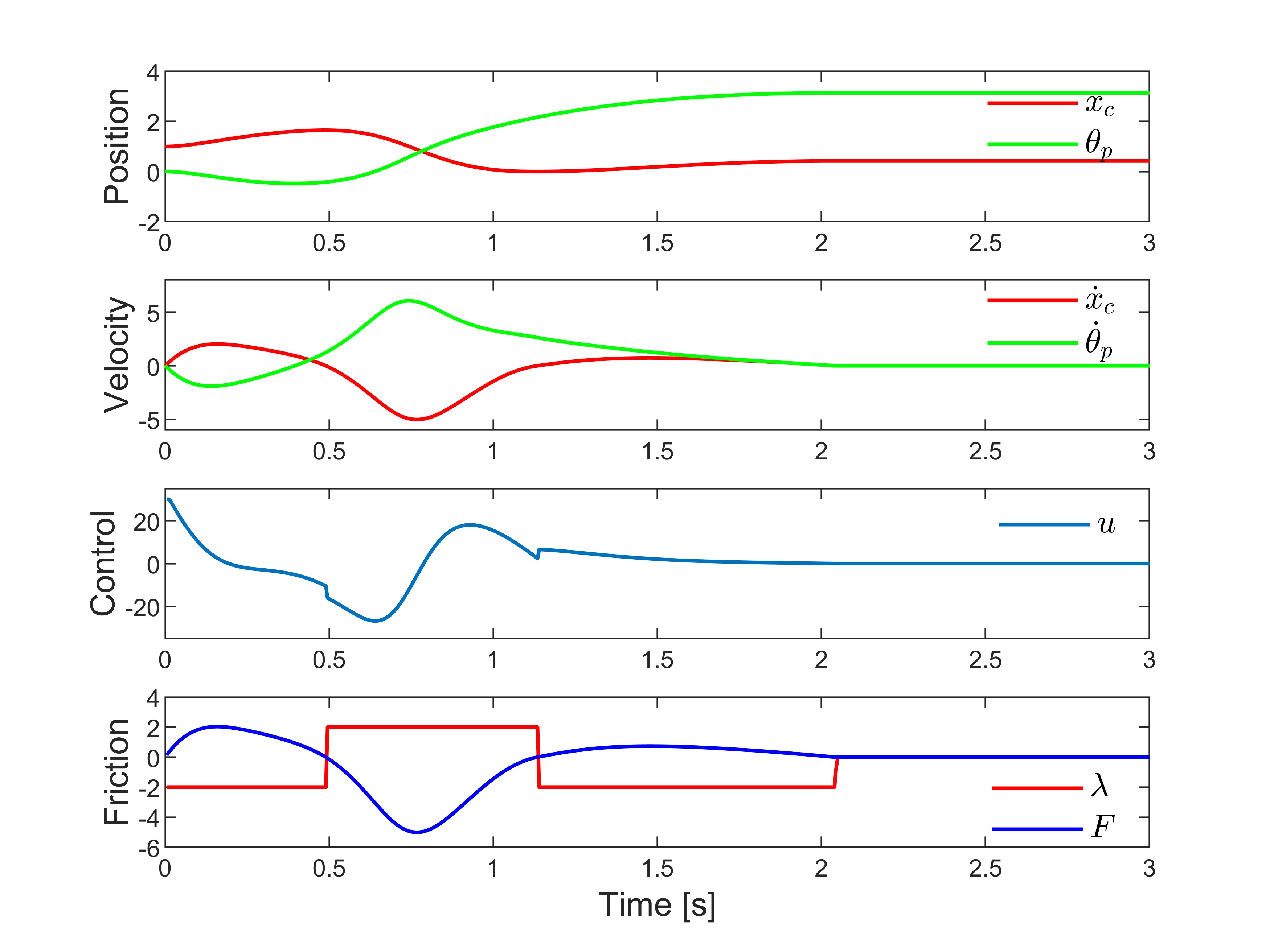}
    \caption{Solution trajectory}
    \label{fig: solution trajectory}
\end{figure}

A solution trajectory obtained by NIP with two Newton corrector steps and $\Delta t = 5 \cdot 10^{-3}$ is depicted in Fig. \ref{fig: solution trajectory}.
The Coulomb friction introduces mode switching between stick and sliding behaviors.
As shown in Fig. \ref{fig: solution trajectory}, the solution trajectory includes four modes without any predefined mode sequence:
\begin{itemize}
    \item From $t = 0$ to $t = 0.490$, the cart slides to the right side to prepare for accelerating;
    \item From $t = 0.490$ to $t = 1.140$, the cart accelerates and slides to the left side to swing the pole;
    \item From $t = 1.140$ to $t = 2.040$, the cart slightly slides to the right side to keep the pole balance;
    \item From $t = 2.040$ to $t = 3$, the cart sticks.
\end{itemize}

\subsection{Necessity of the continuation method}\label{subsection: necessity of the continuation method}
\begin{figure}[!tbp]
    \centering
    \includegraphics[width=0.9\linewidth]{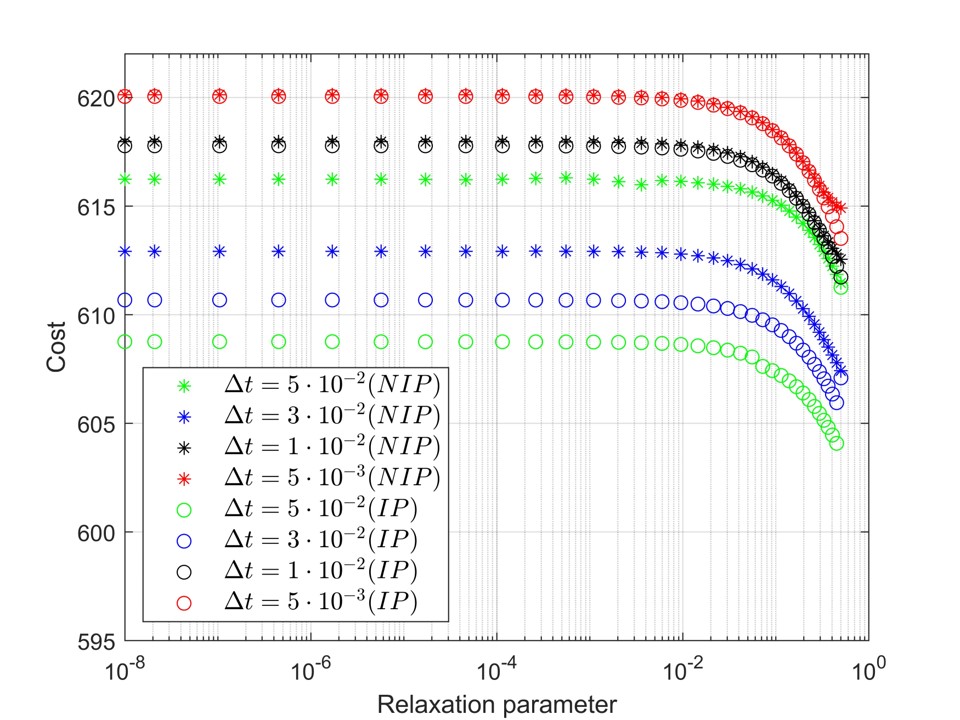}
    \caption{Comparison of different time steps and methods}
    \label{fig: time step comparison}
\end{figure}

In this subsection, we investigate the optimal cost of the discretized OCPEC (\ref{equation: discretized OCPEC}) using different time steps $\Delta t$ and solution methods.
It is well-known that the smaller the time step $\Delta t$, the closer the solution trajectory of the discretized OCPEC (\ref{equation: discretized OCPEC}) approaches the solution trajectory of the continuous-time OCPEC (\ref{equation: continuous OCPEC}).
Since NIP and IP handle the inequality constraint in an entirely different manner, the solution of the first problem $\mathcal{P}(s^0)$ obtained by NIP and IP is generally different, and we cannot expect these two methods to start the solution tracking from the same point.
Nonetheless, as shown in Fig. \ref{fig: time step comparison}, all tests indicate that if the condition $\Delta t \ll s$ holds, then the continuation method can drive the iterates to converge to the optimal solution of the continuous-time OCPEC (\ref{equation: continuous OCPEC}).
However, once the condition $\Delta t \ll s$ fails to hold, the iterates converge to the nearest spurious local solution.
They then remain stuck and unable to escape, regardless of the decreasing relaxation parameter $s$ or the solution method utilized.
In other words, the discretized OCPEC (\ref{equation: discretized OCPEC}) does not approximate the continuous-time OCPEC (\ref{equation: continuous OCPEC}) with correct sensitivity as long as the condition $\Delta t \ll s$ fails to hold.
This observation in the DVI system is identical to the conclusion in the nonsmooth ODE \cite{stewart2010optimal} and dynamical complementarity system \cite{nurkanovic2020limits}, confirming the importance of the condition $\Delta t \ll s$ and the necessity of the continuation method for solving OCPEC. 

\subsection{Performance comparison between NIP and IP}
\begin{figure}[!tbp]
    \centering
    \includegraphics[width=0.9\linewidth]{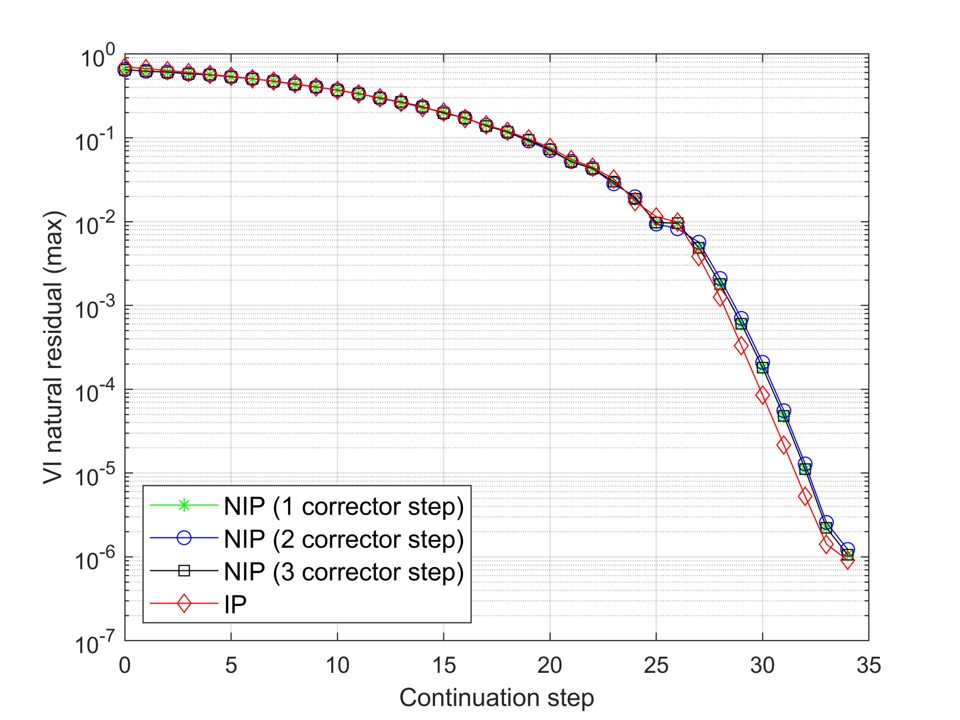}
    \caption{Comparison of VI natural residual}
    \label{fig: VI natural residual comparison}
\end{figure}

\begin{figure}[!tbp]
    \centering
    \includegraphics[width=0.9\linewidth]{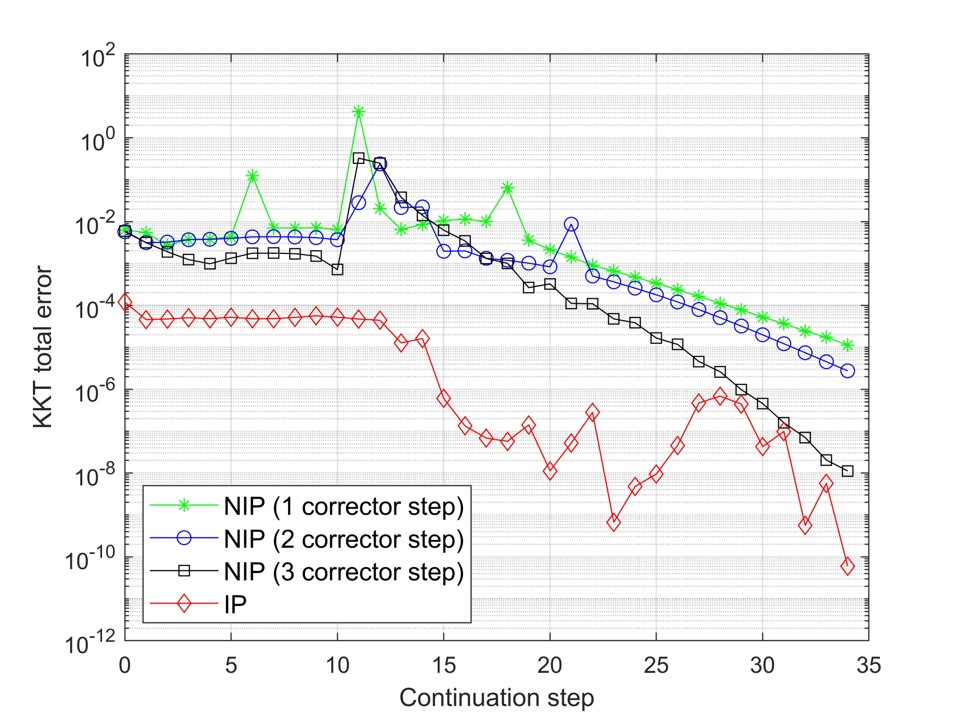}
    \caption{Comparison of KKT error}
    \label{fig: KKT error comparison}
\end{figure}

\begin{figure}[!tbp]
    \centering
    \includegraphics[width=0.9\linewidth]{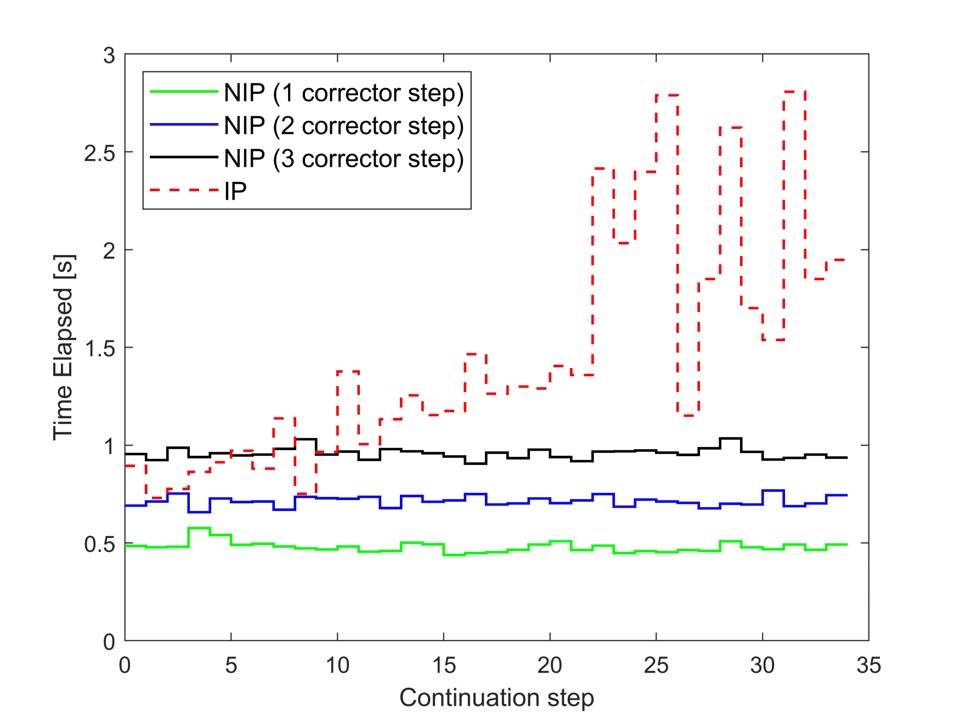}
    \caption{Comparison of elapsed time}
    \label{fig: time Elapsed comparison}
\end{figure}

Since the continuation method is necessary, tracking the solution trajectory accurately and efficiently is of utmost importance.
As shown in Fig. \ref{fig: time step comparison}, the cost trajectory of NIP and IP are nearly identical when $\Delta t = 10^{-2}$; therefore, we compare the performance of NIP and IP using the discretized OCPEC (\ref{equation: discretized OCPEC}) with $\Delta t = 10^{-2}$.
Furthermore, we implement NIP with a different number of Newton corrector steps to compare the solution quality.
As shown in Fig. \ref{fig: VI natural residual comparison}, the violation of equilibrium constraints decreases as the continuation step increases (i.e., $s \rightarrow 0$), and all the implementations track the solution trajectory with a nearly identical violation.
As shown in Fig. \ref{fig: KKT error comparison}, the IP exhibits the smallest KKT error as it solves each problem $\mathcal{P}(s)$ exactly, whereas the KKT error $E^k_{kkt}$ (defined by (\ref{equation: termination condition composite KKT}) and generally dominated by dual residual) of the NIP is larger than that of IP, but can generally be decreased with additional Newton corrector steps.
Due to the nonlinearity and nonsmoothness of the benchmark problem, the Euler predictor step occasionally fails to provide a `good' prediction, thereby hindering NIP in achieving higher-precision solution tracking.
The computation time is the next focal point of attention.
As shown in Fig. \ref{fig: time Elapsed comparison}, the computation time of NIP increases with the number of corrector steps, and the time required for each NIP implementation remains relatively constant in each continuation step.
In contrast, IP demands more time in each continuation step, especially when $s$ is sufficiently small.
Considering all factors into account, NIP with one Newton corrector step emerges as the best choice;
it accurately tracks the solution trajectory and demands significantly less computation time compared to IP.

%%%%%%%%%%%%%%%%%%%%%%%%%%%%%%%%%%%%%%%%%%%%%%%%%%%%%%%%%%%%%%%%%%%%%%%%%%%%%%%%
\section{Conclusion}\label{section: conclusion}
In this study, we proposed a two-stage solution method to solve OCPEC efficiently.
The proposed method relaxes equilibrium constraints to mitigate the numerical difficulties and formulates the KKT system as a perturbed equation system.
In the first stage, a non-interior-point method was employed to estimate an initial point of the solution trajectory.
In the second stage, a predictor-corrector continuation method was utilized to track the solution trajectory. 
The convergence was analyzed.
Numerical experiments revealed that compared with the interior-point-based solvers, the proposed method can accurately track the solution trajectory with significantly less computational time.
Regarding future research, we are considering the embedded real-time application of the proposed method, and attempting to improve the Euler predictor step.

\bibliographystyle{plain} 
\bibliography{reference}

\end{document}